\documentstyle{amsppt}          
\TagsOnRight
\magnification=1100

\define\g{{\goth g}}
\define\h{{\goth h}}

\define\m{{\goth m}}


\define\={\overset\text{def}\to=}
\redefine\B{{\Cal B}}
\define\C{{\Bbb C}}
\define\R{{\Bbb R}}
\define\Z{{\Bbb Z}}

\redefine\D{{\Cal D}}

\redefine\Im{\text{\bf Im}}
\define\Ann{\text{\bf Ann}}
\redefine\ker{\operatorname{Ker }}
\define\GL{{\operatorname{GL}}}
\redefine\P{P_{\text{CM}}(M)}
\define\su{\goth{su}}
\define\Q{{\Cal Q}}

\define\GQ{G_{\Cal Q}}
\define\HQ{H_{\Cal Q}}
\define\gQ{\g_{\Cal Q}}
\define\hQ{\h_{\Cal Q}}
\define\aut{\operatorname{Aut}}
\define\gsym{g_{\operatorname{sym}}}

\TagsOnRight
\topmatter

\title
Explicit construction of a Chern-Moser connection \\for CR manifolds of
codimension two
\endtitle

\author
Gerd Schmalz and Andrea Spiro\\
\phantom{aa}
\endauthor






\leftheadtext{G. Schmalz and A. Spiro}
\rightheadtext{Chern-Moser connection
for CR manifolds of codimension two}

\abstract  In the present paper we suggest an explicit construction of a Cartan
connection for
an elliptic or hyperbolic CR manifold $M$ of dimension six and codimension
two, i.e. a pair $(P, \omega)$, consisting of a principal bundle $\pi: P
\to M$
over $M$ and of a Cartan connection form $\omega$ over $P$, satisfying the
following property:
the (local) CR transformations $f: M\to M$ are in one to one correspondence
with
the (local) automorphisms $\hat f: P \to P$ for which $\hat f^* \omega = \omega$.
For any $x\in M$, this construction determines an explicit monomorphism of the
stability
subalgebra $Lie(\aut(M)_x)$ into the Lie algebra $\h = Lie(H)$ of the
structure
group $H$ of $P$.
\endabstract


\endtopmatter
\document
\subhead 1. Introduction
\endsubhead
\bigskip
A {\it CR structure\/} on a smooth manifold $M$ is a pair $(\D, J)$
where $\D$ is a distribution in $TM$ and $J$ is a smooth family of 
complex structures $J_x: \D_x \to \D_x$ on the spaces $\D_x \subset T_xM$. 
Such geometric structure occurs naturally in studying the geometry 
of 
embedded   submanifolds  of $\C^n$. In fact, typical examples of CR structures are given 
by the pairs $(\D, J)$ on a (sufficiently regular) real submanifold $M \subset \C^n$, where 
$\D$ is the distribution of tangent subspaces $\D_x = \{\ v\in T_xM\ : \sqrt{-1} 
v \in T_xM\ \}$ and  $J$ is the  family of 
 complex structures $J_x: \D_x \to \D_x$ defined by $J_x(v) = \sqrt{-1} v$.
\par
\smallskip
Any CR structure is  naturally associated with the pair of 
distributions 
$\D^{10}, \D^{01}$ in the complexified tangent space $T^\C M$, determined 
by the subspaces $\D^{10}_x$, $\D^{01}_x \subset \D^\C \subset T^\C_xM$ defined by 
$$\D^{10}_x = \{\ v\in \D_x^\C\ :\ Jv = \sqrt{-1} v\ \}\ ,\ 
\D^{01}_x = \{\ v\in \D^\C_x\ :\ Jv = - \sqrt{-1} v\ \}\ .$$
A CR structure $(\D, J)$ is called {\it integrable\/} if the corresponding distribution 
$\D^{10}$ is involutive, i.e. if for any two complex vector fields $X, Y$ with values in 
$\D^{10}$ the Lie bracket $[X,Y]$ is still a complex vector field with values in $\D^{10}$.
The {\it codimension\/} of a CR structure $(\D,J)$ is the codimension 
of the distribution $\D$ in $TM$.\par
A CR structure $(\D,J)$ of codimension one
is called {\it Levi non-degenerate\/} if  the distribution $\D$  is  a contact distribution, 
i.e. if for any 1-form $\theta$ with $\ker \theta_x = \D_x$, $x\in M$, one has that 
$\theta \wedge (d\theta)^n \neq 0$.  A CR structure $(\D, J)$ of arbitrary 
codimension is called {\it Levi non-degenerate\/}   if $\D$ satisfies a 
certain set of conditions, which 
generalizes the condition of Levi non-degeneracy 
of the CR structures of  codimension one. For the exact definition of Levi non-degeneracy 
in arbitrary codimension, see Definition 2.3 below. \par
{\it From now on, any CR structure will be assumed to be integrable and Levi non-degenerate\/}.\par
\smallskip
Now, for a given CR manifold $M$ of codimension $k$, 
we call  {\it osculating 
quadric at a point $x \in M$\/}  the homogeneous real submanifold of $\C^{\frac{\dim M + k}{2}}$
$$\Q_x = \{\ (z,w^1, \dots, w^k) \in \C^{\frac{\dim M + k}{2}}\ :\quad
 \operatorname{Im}(w^i) = \underset{(x)}\to{H}^i(z, \bar z)\ ,\quad 1 \leq i \leq k\ \}\ ,$$
where the $\underset{(x)}\to{H}^i$'s are 
the components of the Levi form of the CR structure, evaluated at the point $x$
(for the definition of Levi form of a CR structure, see \S 2 below). 
Geometrically speaking, the osculating 
quadric at a point $x$ can be characterized as the homogeneous CR manifold, whose CR structure
``osculates up to the second order"
the CR structure of $M$ at the point $x$. We say that a CR manifold $M$ is of {\it strongly uniform type\/}
if all  osculating quadrics $\Q_x$, $x\in M$, are   equivalent. \par
For CR structures of codimension one strong uniformity is an
automatic consequence of non-degeneracy, due to the classification of quadrics
by their signature. In higher codimensions there are only a few cases
where non-degenerate quadrics admit a discrete classification. Among
them most interesting are certainly those with large automorphism groups.
According to \cite{10}, this happens only in two situations: either when
the real dimension of $\Q$
is $6$ and the codimension is $2$, or when the real dimension of $\Q$ is
$2n+n^2$ and the codimension is $n^2$, where $n$ is the complex dimension of 
$\Cal D$. \par
In this paper we discuss the geometry of strongly uniform
CR manifolds of dimension six and codimension two. \par
\medskip
Strongly uniform
CR manifolds  of  dimension six and codimension two have been intensively studied
(see e.g. \cite{14}, \cite{8}, \cite{9}, \cite{10}, \cite{7}, \cite{11}, 
\cite{12}, \cite{16}, \cite{18}, \cite{19}, \cite{6}).
In particular, it is known that for such a  kind of CR manifolds, there are only three possibilities for 
the osculating quadric. Using this fact,  these manifolds are subdivided into three disjoint classes, 
namely  the 
 {\it elliptic\/}, 
{\it parabolic\/} and {\it hyperbolic\/} manifolds. Furthermore, in a previous paper
J. Slov\'ak  and the first author  obtained the following theorem, which is consequence of
the results  in \cite{20} and \cite{5}.
\medskip
\proclaim{Theorem 1.1} \cite{18} Let $M$ be an elliptic or hyperbolic manifold and  $\Q$
 the osculating quadric at one of its points. Denote also by $\GQ$ the  automorphisms group
$\GQ = \aut(\Q)$, by
$\HQ$ the stability subgroup  $\HQ = \aut(\Q)_0$ at $0 \in \Q$,
  and by $\gQ$ the Lie algebra $\gQ = Lie(\GQ)$. 
Then:
\roster
\item"i)" there exists a principal bundle $\pi: P(M) \to M$ with structure group $\HQ$
and a natural injective homomorphism  
$$\imath: \aut_{loc}(M) \to \aut_{loc}(P(M))$$
(called {\rm lifting map\/})
 from 
the pseudogroup $\aut_{loc}(M)$ of local CR transformations of $M$ 
into the pseudogroup $\aut_{loc}(P(M))$ of local automorphisms of 
bundle $P(M)$, which makes the following diagram commuting
$$
\CD
 P(M) @>\imath(f)>>
P(M)  \\
@V\pi VV @VV\pi V\\
M@>f>> 
M
\endCD
$$
\item"ii)" there exists  
a Cartan connection
$\omega_M: TP(M) \to \gQ$  such that 
a diffeomorphism $\tilde f: P(M) \to P(M)$ is equal to $\tilde f = \imath(f)$  for
some $f\in \aut(M)$
if and only if $\tilde f^*\omega_M = \omega_M$. 
\endroster
\endproclaim
\medskip
Notice  that, by Theorem 2.7 in \cite{20}, if  $(P'(M), \omega'_M)$ is a pair 
consisting of a principal bundle $P'(M)$, which  satisfies (i),
 and of a  Cartan connection $\omega'_M$, which  
satisfies (ii)  plus a certain set of linear conditions on the  curvature, 
then $(P'(M), \omega'_M)$ is isomorphic with the pair $(P(M), \omega_M)$ of 
Theorem 1.1 (i.e. there exists a principal bundle isomorphism $f: P'(M) \to P(M)$
such that $f^* \omega_M = \omega'_M$). \par
\medskip
From the existence of the invariant 
Cartan connection $\omega_M$, 
it follows that  the action of the stability subgroup 
$\aut(M)_x$, $x\in M$,  on the fiber $P_x = \pi^{-1}(x) \subset P(M)$,  commutes with the 
simply transitive action of $\HQ$ on $P_x$. From this it may be inferred   
that, given a point $u\in P_x$, 
the map 
$$X\in Lie(\aut(M))_x \mapsto V_{u,X}\in Lie(\HQ)\ ,$$
where $V_{u,X} \in Lie(\HQ)$ is  the element, whose 
 fundamental vector field $V^*_{u,X}$
satisfies $V^*_{u,X}|_u = \imath_*(X)|_u$, 
is a Lie algebra injective homomorphism from  $Lie(\aut(M)_x)$ into 
$Lie(\HQ)$. 
By this remark,  the explicit knowledge of the bundle  $P(M)$ and of the lifting map $\imath$
gives an effective tool for reducing several questions 
on elliptic or 
hyperbolic manifolds to questions concerning some special subalgebras of $Lie(\HQ)$.\par 
\smallskip
It is a 
disadvantage of the iterative construction, on which the result in
\cite{18} is
 based,  that it does not provide an explicit expression of the bundle $P(M)$
and  of the Cartan connection (for more details, see also \cite{5}, \cite{18}
or the review in \cite{1}). On the other hand, 
the explicit construction of bundle and absolute parallelism  in \cite{7}
does not give, in general,  a Cartan connection and hence, it does not give an explicit 
expression for the homomorphism
between $Lie(\aut(M)_x)$ and $Lie(\HQ)$.\par
\medskip
In this paper we give an explicit and (in our opinion) simple
construction of a  Cartan connection $(\P, \psi_{CM})$,  which satisfies all the 
claims of Theorem 1.1.
However, a direct check shows  
 that, generically,  the curvature of $\psi_{CM}$ 
{\it does not\/} satisfy the linear equations required by 
Theorem 2.7 in  \cite{20}.  This means that, generically, 
the new pair $(\P, \psi_{CM})$ {\it is not
isomorphic with the pair $(P_M, \omega_M)$\/}. However, we expect that 
our bundle $\P$ admits also another  
Cartan connection, which  does 
 satisfy the conditions 
of  the quoted theorem. If this conjecture is correct, this would
imply that the bundle $\P$ is  equivalent 
to the bundle $P(M)$. 
\par
\medskip
The bundle $\P$ is  the exact analogue of the principal
 bundle associated with a Levi non-degenerate CR manifold of codimension one, 
introduced by S.S. Chern and J. Moser in \cite{4}. 
The construction can be roughly described as follows. 
First of all, we consider a principal bundle $\pi_1: E(M) \to M$, with 2-dimensional fibers, 
formed by all pairs of 1-forms $(\theta^1, \theta^2) \in T^*_xM \times T^*_xM$, 
$x\in M$, such that $\D_x = \ker \theta^1 \cap \ker \theta^2$ and so that 
some additional condition, which are specified below,  on the corresponding Levi form 
are satisfied. Secondly, we consider a special class of 
linear frames $u = (e_1, \dots, e_8)$ at the tangent spaces of $E(M)$, which we call
{\it adapted to the CR structure\/}. The 
bundle $\pi_2: \P \to E(M)$ of adapted frames is proved to 
have a natural structure of  an $\HQ$-principal bundle  $\pi = \pi_1 \circ \pi_2: \P \to M$
 over $M$ and it satisfies claim (i) of Theorem 1.1.
In particular, since any CR diffeomorphism $f: M \to M$
lifts  naturally to a map $\hat f: E(M) \to E(M)$ and  each map $\hat f$ lifts 
naturally to a 
diffeomorphism $\Hat{\Hat f}: L(E(M)) \to L(E(M))$ of the linear frame bundle 
$L(E(M))$, the  lifting map $\imath: \aut(M) \to \aut(\P)$ of
 Theorem 1.1 (i)
is simply the map $\imath: f \mapsto \Hat{\Hat f}$.\par
Even the construction of the Cartan connection $\psi_{CM}$ is
modelled on the arguments used in \cite{4}. In fact, 
starting from a set of 1-forms $\{\omega^i, \psi^j\}$,
which represents an arbitrary $\gQ$-valued Cartan connection on $\P$, we show 
how to determine  a set of real functions $S_i^j$ on $P(M)$
 so that the new collection of 1-forms
$\{ \omega^i, \psi^j - \sum_i S_i^j \omega^i\}$ represents  a Cartan connection $\psi_{CM}$,
which verifies Theorem 1.1 (ii).  
\par
\medskip
We have to recall that, in case $M$ is real analytic and it is presented 
as a real submanifold of $\C^4$, another explicit embedding
of $Lie(\aut(M)_x)$ into 
$Lie(\HQ)$ can be  obtained  by studying the normal form of the defining equations for 
$M$ (see  \cite{14}, \cite{9}, \cite{12}). On the other hand, the approach given here 
is valid for any smooth manifold $M$ and it is based only on the intrinsic CR geometry of $M$. \par 
\medskip
Before concluding, we would like to point out that the construction of the Chern-Moser bundle 
can be done also in case $M$ is a parabolic manifold. But, in this case, 
the bundle  does not have a natural structure of a principal bundle over $M$. This 
is certainly consistent with the results of \cite{20} and \cite{5}, which 
imply  neither   existence  
nor non-existence of a  canonical Cartan connection 
on parabolic manifolds. At the best of our knowledge, 
it is not known if there is  any 
obstruction to the existence of a canonical Cartan connection on parabolic manifolds.\par
\medskip
The plan of the paper is the following. In \S 2, we discuss a few general facts on CR structures
and  we give the definition of 
strongly uniform CR manifolds of dimension six and codimension two. We also give a new proof of the 
fact that 
these manifolds are subdivided in 
exactly three classes. The approach we use  is strictly related (even if  
independent) with the discussion of hyperbolic and elliptic distributions of codimension 2
by  A. \v Cap and 
M. Eastwood  in \cite{3}. \par
In \S 3 and \S 4, we construct the  bundle $\P$ of an arbitrary elliptic  or 
hyperbolic  manifold $M$. We call it {\it Chern-Moser bundle of $M$\/}. 
In \S 5 we prove that $\P$ admits a natural structure of 
principal bundle over $M$ and we show that its structure group 
is $\GQ$. In \S 6, we construct the  Cartan connection $\psi_{CM}$ on $\P$ and 
we illustrate why the pair
  $(\P, \psi_{CM})$ is generically not isomorphic 
with  the pair $(P(M), \omega_M)$ 
of Theorem 1.1. \par
\smallskip
For any real vector space $V$, we denote by $\GL(V,\R)$  the group of  
linear isomorphisms of $V$ into itself. If there is a complex structure $J$ on $V$, 
we denote by $\GL(V,\C) \subset \GL(V,\R)$  the subgroup of isomorphisms  commuting 
with $J$.\par
\bigskip
\bigskip
\subhead  2. First definitions and preliminaries
\endsubhead
\bigskip
\subsubhead 2.1 Distribution of contact type 
and Levi non-degenerate CR structures of codimension $k$\endsubsubhead
\par
\medskip
Let $M$ be a manifold of dimension $n$ and let $\D\subset TM$ be a distribution
of codimension $k$ on $M$. For any  point $x\in M$, 
we call {\it conormal frame at $x$\/} any k-tuple
$\theta = (\theta^1, \dots, \theta^k)$ of linearly 
independent 1-forms $\theta^i\in T^*_xM$ such that
$\bigcap_{i = 1}^k \ker \theta^i =  \D_x$. It is clear that any conormal frame is a basis 
for the subspace of 1-forms vanishing on $\D_x$.  We call 
{\it conormal frame bundle of $\D$\/} the bundle $\pi: E(M, \D) \to M$ of all conormal frames at the 
points of  $M$ (here,  $\pi$ is the natural 
projection map which sends any conormal frame  $\theta \in T^*_xM \times \dots \times T^*_xM$ 
to the point $x$). 
Notice that the right action of $\GL_k(\R)$  on $E(M, \D)$ defined by 
$$\mu: \GL_k(\R)\times E(M, \D) \to E(M, \D)\ ,\ \mu(A, (\theta^1, \dots, \theta^k)) \= 
((A^{-1})^1_j \theta^j, \dots, (A^{-1})^k_j \theta^j)$$
acts transitively on the fibers of $\pi: E(M, \D)\to M$ and makes 
$E(M, \D)$  a principal bundle over $M$. \par
\medskip
With any conormal frame $\theta = (\theta^1, \dots, \theta^k) \in E(M, \D)_x$, $x\in M$, we may associate a k-tuple 
of 2-forms in $\Lambda^2\D_x$ as follows: consider a smooth section $\tilde \theta: \Cal U \subset M \to E(M, \D)$
defined on a neighborhood $\Cal U$ of $x$, such that $\tilde \theta_x = \theta$; then let us define
$$\widetilde{d\theta} = (\widetilde{d\theta}^1, \dots, \widetilde{d\theta}^k) \in \Lambda^2 \D_x \times \dots \times 
\Lambda^2 \D_x\ ,$$
$$\widetilde{d\theta}^i(X,Y) \= d(\tilde \theta^i)_x(X,Y)\ ,\qquad \text{for any}\ X,Y\in \D_x\ .$$
It can be checked that the 2-forms $\widetilde{d\theta}^i \in \Lambda^2 \D_x $ are independent of the 
choice of the extension $\tilde \theta$ and  that they depend uniquely on  $\theta = \tilde \theta_x$. 
We will call $\widetilde{d \theta}$ the {\it $\R^k$-valued 2-form associated with \/} $\theta$.\par
\medskip 
\definition{Definition 2.1} 
A codimension $k$ distribution $\D \subset TM$ 
is called  {\it of contact type\/} if for any point $x\in M$  and for any 
conormal frames $\theta = (\theta^1, \theta^2, \dots, \theta^k) \in E(M, D)$
\roster
\item"a)" $\Im(\widetilde{d\theta}) = \operatorname{Span}
\{\ v\in \R^k\ :\ v = \widetilde{d\theta}(X,Y)\ \text{for some}\ X,Y\in \D_x\ \} = \R^k$;
\item"b)" $\Ann(\widetilde{d\theta}) =  \{\ X \in \D_x\ :\ \widetilde{d\theta}(X,*) = 0\ \} = \{0\}$.
\endroster
\enddefinition
\bigskip
A {\it CR structure of codimension k on a manifold $M$\/}
 is a pair $(\D, J)$, where 
$\D$ is a distribution of codimension $k$ and $J$ is  a smooth family of
complex structures
$$J_x : \D_x \to \D_x$$
which satisfy the integrability conditions:
$$J([JX, Y] + [X, JY])\in \D\ ,\tag{$2.1_1$} $$
$$[JX, JY] - [X, Y] - J([JX, Y] + [X, JY])  = 0\ .\tag{$2.1_2$}$$
We recall that the integrability conditions (2.1)   
are satisfied if and only if  the eigendistributions
$\D^{10} \subset T^\C M$ and $\D^{01} \subset T^\C M$ of $J$, given 
by the $J$-eigenspaces in $\D^\C$ corresponding with the eigenvalues $i$ and $-i$,
are involutive, i.e.  the space of their sections is closed
under Lie brackets.\par
If $(\D, J)$ is a CR structure of codimension $k$ and 
$\theta = (\theta^1, \dots, \theta^k)$ is a conormal frame at a point $x\in M$, 
we call
{\it Levi form of $(\D, J)$ associated with $\theta$\/} the k-tuple of 
bilinear forms on $\D_x$ defined by 
$$L^\theta = (L^{\theta^1}, \dots, L^{\theta^k}) \ ,\qquad 
L^{\theta^i}(X, Y) \= \widetilde{d\theta^i}(X, JY)\ ,\quad \text{for any}\ X, Y\in \D_x\ .\tag 2.2$$
Using the integrability conditions (2.1), one can check that for any $i$ and for any $X, Y\in \D_x$
$$L^{\theta^i}(JX,JY) = L^{\theta^i}(X,Y)\ ,\qquad L^{\theta^i}(X,Y) = L^{\theta^i}(Y,X)\ ,\tag 2.3$$
and hence  $L^\theta$ is an $\R^k$-valued, $J$-invariant, symmetric, bilinear form on $\D_x$.\par
We also call {\it complex Levi form at $x$\/} the  
$\C^k$-valued 
hermitian form defined by 
$$\Bbb L^\theta(Z, W) \= \frac{1}{2i} d\tilde \theta_x(Z, \overline{W})\ ,\tag 2.4$$
for any $Z, W \in \D^{\C}_x$. It is quite simple to check that  $\Bbb L^\theta$ 
coincides (up to some isomorphism between $\C^k$ and 
$T^\C_x/D^\C_x$) with the  complex  Levi form at $x$ as classically defined (see e.g. \cite{2})
and that  the 
 Levi form $L^\theta$ is (up to a factor)  the real  part 
of the  complex Levi form. 
\par
\medskip
\definition{Definition 2.2} A CR structure $(\D,J)$ of codimension $k$ is called
{\it Levi non-de\-ge\-ne\-ra\-te\/}\ if the underlying distribution $\D$ is 
of contact type.
\enddefinition
\medskip
Notice that from definitions and (2.4), a CR structure 
is Levi non-degenerate if and only if any 
$\C^k$-valued Levi form $\Bbb L^\theta$ 
satisfies 
$$\Im (\Bbb L^\theta_x) = 
\operatorname{Span}\{\ v\in \C^k :\ v = \Bbb L^\theta_x(Z,W)\ \text{for some}\ Z,W\in \D^\C_x\ \} = \C^k\ ,\tag
2.5$$
$$\Ann (\Bbb L^\theta_x )\cap \D^{10}_x = 
\{\ Z \in \D^\C_x\ :\ \Bbb L^\theta_x(Z, * ) = 0\ \} \cap \D^{10}_x = \{0\}\ .
\tag 2.6$$
In particular,  it is clear that  Definition 2.2 is simply a reformulation of the
usual  Levi non-degeneracy condition on complex valued Levi forms (see e.g. \cite{7}).
\par
\bigskip
\subsubhead 2.2.  The three types of CR manifolds of dimension six and  codimension two
\endsubsubhead
\bigskip
From now, by $(M, \D, J)$ we will always denote 
an integrable Levi non-de\-ge\-ne\-ra\-te CR manifold of dimension six and
codimension two. \par
Consider such a  CR manifold $(M,\D, J)$ and 
let $\pi: E(M, \D) \to M$ be the associated bundle of conormal frames of the  distribution $\D$. 
Recall that $E(M,\D)$ is a $\GL_2(\R)$-principal bundle over $M$.
The aim of this subsection is to construct a $\GL_2(\R)$-equivariant map 
$$\Psi: E(M,\D) \to \Bbb P(S_{2\times 2}(\R))\ ,$$
where $\Bbb P(S_{2\times 2}(\R))$ denotes the projective space of the space 
of symmetric $2 \times 2$ real matrices 
and on which we consider the right action of $\GL_2(\R)$ given by
$$A\cdot [a] = [ (A^{-1})  \cdot a\cdot (A^{-1})^T]\qquad \text{for any}\ A \in \GL_2(\R)\ .\tag 2.7$$
Since  
the map $\Psi$ will be proved to be $\GL_2(\R)$-equivariant and any fiber $E_x = \pi^{-1}(x) \subset E(M,\D)$
is a $GL_2(\R)$-orbit in $E(M,\D)$, we will have a well  defined map which assigns to
any point $x\in M$ the  $\GL_2(\R)$-orbit in $\Bbb P(S_{2\times 2}(\R))$ containing
the image $\Psi(E_x)$ of  $E_x = \pi^{-1}(x)$.\par
\smallskip
In order to define the map $\Psi$, we first need to recall
the following   lemma. \par
\medskip
\proclaim{Lemma 2.3} Let $(V, J)$ be a 4-dimensional real vector space endowed with a 
complex structure $J$ and let $\Lambda^2_\C V^*$ be the space of $J$-invariant 
2-forms of $V$. Let also $\tau \in \Lambda^4 V$ be a non-trivial 4-vector on $V$ and 
${\underset{(\tau)}\to G}$ the symmetric bilinear form defined by 
$${\underset{(\tau)}\to G}: \Lambda^2_\C V^* \times \Lambda^2_\C V^* \to \R
\ ,\qquad {\underset{(\tau)}\to G}(\alpha, \beta) \= (\alpha\wedge \beta)(\tau)\ .\tag 2.8$$
Then ${\underset{(\tau)}\to G}$ is a Lorentz metric on $\Lambda^2_\C V^*$ and the image 
$\rho(\GL_2(V,\C))$ of the representation 
$$\rho: \GL_2(V,\C) \to \GL(\Lambda^2_\C V^*, \R)\ ,\qquad 
\rho(A) \cdot \alpha = A^* \alpha\ ,\tag 2.9$$
is equal to the connected component of the identity of 
the linear conformal group of $(\Lambda^2_\C V, {\underset{(\tau)}\to G})$.
\endproclaim
\demo{Proof} Consider a  basis $(e_1, \dots, e_4)$ for $V$, which satisfies  $J e_1 = e_2$,  $Je_3 = e_4$ 
and $e_1\wedge \dots \wedge e_4 = \tau$, and let
$(e^1, \dots, e^4)$ be the corresponding dual basis in $V^*$. Then we may consider the following basis 
for $\Lambda^2_\C V^*$   
$$\xi^0 = \frac{1}{2}\left(e^1 \wedge e^2 + e^3 \wedge e^4\right)\ ,\qquad 
\xi^1 = \frac{1}{2}\left(e^1 \wedge e^2 - e^3 \wedge e^4\right)\ ,\tag 2.10$$
$$\xi^2 = \frac{1}{2}\left(e^2 \wedge e^4 + e^1 \wedge e^3\right)\ ,\qquad 
\xi^3 = \frac{1}{2}\left(e^1 \wedge e^4 - e^2 \wedge e^3\right)\ .\tag 2.11$$
and observe that the matrix associated with the bilinear form ${\underset{(\tau)}\to G}$ in this basis is
$$\left({\underset{(\tau)}\to G}(\xi^i, \xi^j)\right) = \left( 
\matrix 1 & 0 & 0 & 0 \\
0 & -1 & 0 & 0 \\
0 & 0 & -1 & 0 \\
0 & 0 & 0 & -1
\endmatrix \right)\ .$$
i.e. ${\underset{(\tau)}\to G}$ is a Lorentz metric. 
For any complex endomorphism $A \in \GL(V, \C)$, let us denote
by $A^*$ the induced action on $\Lambda^2_\C V^*$. Then, for any
$\alpha \in \Lambda^2_\C V$,
$${\underset{(\tau)}\to G}(A^*\alpha, A^*\alpha) =
\left|\det(A)\right|^2 {\underset{(\tau)}\to G}(\alpha, \alpha)\ ,$$
where by $\det(A)$ we mean the determinant of the complex matrix associated 
with $A$ in the complex basis  $(e_1 - i J e_1, e_3 - i Je_3)$. 
Then,  for any 
$\alpha \in \Lambda^2_\C V$, we have that $A^*\alpha \in \Lambda^2_\C V$
and that  
$${\underset{(\tau)}\to G}(A^*\alpha, A^*\alpha) =  
\left|\det(A)\right|^2 {\underset{(\tau)}\to G}(\alpha, \alpha)\ .$$
This shows that  (2.9) determines a representation of $\GL(V, \C)$ as a subgroup
of the linear conformal group of ${\underset{(\tau)}\to G}$. Since the kernel of this 
representation is $\text{SO}_2(\R)$, by counting dimension, we get that $\rho(\GL(V, \C))$ is isomorphic
to the connected component of the identity of $\text{CO}_{3,1}(\R)$.\qed
\enddemo
\bigskip
Now, for any given point $x\in M$, let us fix a non trivial element $\varpi \in \Lambda^4 \D_x$
and let us consider the  
Lorentz  metric  $\underset{(\varpi)}\to G$ on $\Lambda^2\D_x$  defined in (2.8). Notice that if we replace 
$\varpi$ by some other non trivial element of $\Lambda^4 \D_x$, the associated Lorentz metric 
changes only by multiplication by a non trivial factor.\par
Let also $E_x \subset E(M, \D)$ be
 the fiber over $x$ in $E(M,\D)$ 
and, for any $\theta = (\theta^1, \theta^2)\in E_x$, denote by $a(\theta) = (a^{ij}(\theta))$ 
the symmetric matrix
$$a^{ij}(\theta)  = {\underset{(\varpi)}\to G}(\widetilde{d\theta}^i, \widetilde{d\theta}^j)\ .\tag 2.12$$
Using the condition of Levi non-degeneracy, one can check that $a(\theta) \neq 0$ for any 
$\theta \in E_x$ and, by the previous remark, we have that 
 the projective class $[a(\theta)]\in \Bbb P(S_{2\times 2}(\R^2))$ is independent of the choice of $\varpi$. 
In particular, we have a 
well-defined  map 
$$\Psi: E(M, \D)  \to \Bbb P(S_{2\times 2}(\R))\ ,\qquad \Psi(\theta) = [a^{ij}(\theta)]\ .\tag 2.13$$
If we consider on $\Bbb P(S_{2\times 2}(\R))$ the right action given in (2.7), 
 it follows from definitions that $\Psi$ is $\GL_2(\R)$-equivariant and that
it is the map we announced at the beginning of this subsection.\par
\bigskip
Observe that the action (2.7) of  $\GL_2(\R)$ on $\Bbb P(S_{2\times 2}(\R))$
has exactly three distinct 
$\GL_2(\R)$-orbits, namely
 the orbits of the following three elements
$$a) \ \left[\matrix 1 & 0\\ 0 & 1 \endmatrix\right]\ ,\qquad
b) \ \left[\matrix 0 & 0\\ 0 & 1 \endmatrix\right]\ ,\qquad
c) \ \left[\matrix 1 & 0\\ 0 & - 1 \endmatrix\right]\ .\tag 2.14$$
This fact leads immediately to the following definition.\par
\bigskip
\definition{Definition 2.4} We say that  $x \in M$ is a point of {\it elliptic\/}, 
{\it parabolic\/} or {\it hyperbolic type\/} if the image  $\Psi(E_x)$ of the fiber
 $E_x = \pi^{-1}(x) \subset E(M, \D)$ is contained in 
the $\GL_2(\R)$-orbit of  $\left[\matrix 1 & 0\\ 0 & 1 \endmatrix\right]$, 
$\left[\matrix 0 & 0\\ 0 & 1 \endmatrix\right]$ or 
$\left[\matrix 1 & 0\\ 0 & - 1 \endmatrix\right]$, respectively.\par
We say that $M$ has  {\it strongly uniform type\/}
if  all points $x\in M$ have the same type, or, equivalently, if the image
 $\Psi(E(M,\D)) \subset \Bbb P(S_{2\times 2}(\R^2))$ 
is contained in exactly one $\GL_2(\R)$-orbit of $\Bbb P(S_{2\times 2}(\R^2))$. In this case 
$M$ is called {\it elliptic\/}, {\it parabolic\/} or {\it hyperbolic manifold\/} according to 
the type of its points.
\enddefinition
\bigskip
Now, if  $x\in M$ is a point of  a $(2n+k)$-dimensional Levi non-degenerate CR manifold
$(M, \D, J)$  of CR codimension k and $\theta = (\theta^1, \dots, \theta^k)$ is a conormal frame 
at the point $x$, we may consider the
$\C^k$-valued Hermitian form 
$$H_x : \C^{n} \times \C^{n} \to \C^k\ , \qquad 
H_x(z,z') \= \frac{1}{2i} \widetilde{d\theta}(z,z')\ ,$$
where we have  identified $\C^n \simeq \D^{10}_x$. This Hermitian form 
determines uniquely (up to isomorphisms) the so called 
{\it osculating quadric of $M$ at the point  $x$\/}, namely 
the quadric 
$$Q_x \= \{ (z,w) \in \C^{n+k}\ : \ \text{Im} w = H_x(z,z)\ \}\ . \tag 2.15$$
Usually (see e.g. \cite{7}), a CR manifold is called {\it strongly uniform\/} if the 
osculating quadrics $Q_x$ and $Q_{x'}$ of any two points are equivalent, i.e. 
if  $H_{x'} = B \cdot (A^*H_x)$ for some $A\in \GL_n(\C)$
and some $B\in \GL_k(\R)$. It turns out that  if $M$ is of dimension 6 and CR codimension 2, 
$M$ is of strongly uniform type (according to Definition 2.4) if and only if it is strongly uniform
according to this last definition. Actually, we have the following fact, whose
proof can be found e.g. in \cite{16}.\par
\bigskip
\proclaim{Proposition 2.5}
 Let $M$ be a   Levi non-degenerate, integrable, 
CR manifold of dimension six and
codimension two of strongly uniform type. Then, at any point $x\in M$, the 
osculating quadric is (up to equivalence)   exactly one of the following three quadrics
$$\matrix
a) & \quad  \operatorname{Im}(w_1) = \operatorname{Im}(z_1 \bar z_2)& ,& \qquad \operatorname{Im}(w_2) = \operatorname{Re}(z_1 \bar z_2)& ;\\
\ &\ & \ &\ & \ \\
b) &\quad  \operatorname{Im}(w_1) = |z_1|^2& , &\qquad \operatorname{Im}(w_2) = \operatorname{Re}(z_1 \bar z_2)& ;\\
\ &\ & \ &\ & \ \\
c)& \quad  \operatorname{Im}(w_1) = |z_1|^2& ,&\qquad \operatorname{Im}(w_2) = |z_2|^2& .
\endmatrix$$
In particular,   a) occurs  if $M$ is elliptic,  
b) occurs  if $M$ is  parabolic and  c) occurs if 
$M$ is hyperbolic.
\endproclaim
\medskip
We want to  stress the fact that, for a CR structure $(\D, J)$,  
 the property  of being elliptic, parabolic 
or hyperbolic is mainly a quality of the underlying distribution $\D$ and it is  
independent of the nature of complex structure $J$.
See \cite{3}\ for some investigations on  this fact.  \par
\bigskip
From this point on, we will always assume that $(M, \D, J)$ is a strongly uniform
CR manifold of dimension six and  codimension two and that it is either elliptic or hyperbolic.
\bigskip
\bigskip
\subhead 3. A reduction of the conormal frame bundle $E(M,\D)$
\endsubhead
\bigskip
Consider a strongly uniform CR manifold $(M, \D, J)$ of elliptic or 
hyperbolic type and let $\pi: E = E(M, \D) \to M$ be the conormal frame 
bundle determined by the distribution $\D$. 
In all the following, we will denote by 
$\varpi = (\varpi^1, 
\varpi^2)$ the {\it tautological pair of $E$\/}, that is the pair of 1-forms $\varpi^1$ and 
$\varpi^2$ defined at any point $\theta = (\theta^1, \theta^2) \in E_x \subset E(M, \D)$
by 
$$\varpi^a|_\theta(X, Y) = \theta^a(\pi_*(X), \pi_*(Y))\ ,\ a=1,2\ ,
\qquad \text{for all}\ X,Y\in T_\theta E\ .\tag 3.1$$  
Then, we have the following: \par
\proclaim{Proposition 3.1}
 Let $\Psi: E  \to \Bbb P(S_{2\times 2}(\R))$
be the map defined in (2.13) and let  $\hat E = \Psi^{-1}\left([a_o]\right)$, where 
$[a_o] = \left[\matrix 1 & 0\\ 0 & 1 \endmatrix\right]$ 
or $\left[\matrix 0 & 1\\ 1 & 0 \endmatrix\right]$ if $M$ is 
elliptic or hyperbolic, respectively. Then:\par
\roster
\item for any point $x\in M$ and any conormal frame $(\theta^1, \theta^2) \in \hat E_x$, 
there exist four 1-forms $(e^3, e^4, e^5, e^6)$ such that 
$(\theta^1, \theta^2, e^3, \dots, e^6)$ is a basis for $T^*_xM$ such that $e^4|_{\D_x} =
J^* e^3|_{\D_x}$, 
$e^6|_{\D_x} = J^* e^5|_{\D_x}$ and  the
 2-forms $(\widetilde{d\theta^1}, \widetilde{d\theta^2})$
can be written in terms of this basis as in the following table: 
\smallskip
\moveleft 0.3 cm
\vbox{\offinterlineskip
\halign {\strut\vrule\hfil\ $#$\ \hfil
 &\vrule\ \vrule\hfil\ $#$\
\hfil&\vrule\hfil\ $#$\
\hfil&
\vrule\ \vrule\hfil\ $#$\
\hfil \vrule\cr
\noalign{\hrule}
\phantom{\frac{\frac{1}{1}}{\frac{1}{1}}} M
\ \ &
\widetilde{d\theta^1}
{}^{\phantom{\frac{\frac{1}{1}}{\frac{1}{1}}}}
_{\phantom{\frac{\frac{1}{1}}{\frac{1}{1}}}}
&
\widetilde{d\theta^2}
{}^{\phantom{\frac{\frac{1}{1}}{\frac{1}{1}}}}
_{\phantom{\frac{\frac{1}{1}}{\frac{1}{1}}}}
&
\phantom{\frac{\frac{1}{1}}{\frac{1}{1}}}
\hat G
\cr \noalign{\hrule}
\text{Elliptic}
&
\left(e^4\wedge e^6 + e^3\wedge e^5\right)_{\D_x}
{}^{\phantom{\frac{\frac{1}{1}}{\frac{1}{1}}}}
_{\phantom{\frac{\frac{1}{1}}{\frac{1}{1}}}}
&
\left(e^3\wedge e^6 - e^4\wedge e^5\right)_{\D_x}
&
\smallmatrix \phantom{a}\\ 
\left\{\ A \in M_{2,2}(\R)\ :\ 
A = \left(\smallmatrix a & b\\ -b & a\endsmallmatrix \right) \cdot E\ ,
\right.\\ 
 \phantom{a}\\
\text{where}
\ E = \left(\smallmatrix 1 & 0\\ 
0 & 1\endsmallmatrix \right)\ \text{or}\ 
 E = \left(\smallmatrix 0 & 1\\ 
1 & 0\endsmallmatrix \right)\ ,\\
\phantom{a}\\
\left.\text{and}\ a^2 + b^2 \neq 0 \right\}\ \simeq \ \C^* \times \Bbb Z_2 \\
\phantom{a}
\endsmallmatrix
\cr \noalign{\hrule}
\text{Hyperbolic}
&
\left(e^3 \wedge e^4\right)_{\D_x}
{}^{\phantom{\frac{\frac{1}{1}}{\frac{1}{1}}}}
_{\phantom{\frac{\frac{1}{1}}{\frac{1}{1}}}}
&
\left(e^5\wedge e^6\right)_{\D_x}
&
\smallmatrix \phantom{a}\\ 
\left\{\ A \in M_{2,2}(\R)\ :\ 
A = \left(\smallmatrix a & 0\\ 
0 & b\endsmallmatrix \right) \cdot E\ ,\right.\\
\phantom{a}\\ 
\text{where}\ E = \left(\smallmatrix \pm 1 & 0\\ 
0 & \pm 1\endsmallmatrix \right)\ \text{or}\ 
 E = \left(\smallmatrix 0 & 1\\ 
1 & 0\endsmallmatrix \right)\\
\phantom{a}\\
\left.\text{and}\ a, b\in \R_{\geq 0}\ \right\}\ \simeq\ \R_* \times \R_* \times \Z_2\\
\phantom{a}
\endsmallmatrix
\cr \noalign{\hrule}
}}
\smallskip
\centerline{\bf Table 1}
\smallskip
\item  $\hat E \subset E$ is a reduction of $\pi: E \to M$ with 
structure group $\hat G \subset \GL_2(\R)$ given 
 in Table 1;
\item let $\tilde G\subset \GL_2(\D_x)$ and 
$\rho: \tilde G \to \hat G$ be the subgroups and isomorphisms given in 
Table 2; then, for any element $\theta \in \hat E_x \subset E$ and for any $A \in \tilde G$,
we have that $\widetilde{d\theta} \circ \tilde A = \rho(A) \cdot \widetilde{d\theta}$
(in the following table, the elements of $\tilde G$ are identified with 
the matrices which are associated to their action on $\D^*_x$ w.r.t. a basis 
$e^3, e^4, e^5$ and $e^6$ satisfying (1)):
\par
\medskip
\moveleft 0.6 cm
\vbox{\offinterlineskip
\halign {\strut\vrule\hfil\ $#$\
\hfil&\vrule\hfil\ $#$\
\hfil&
\vrule\ \vrule\hfil\ $#$\
\hfil \vrule\cr
\noalign{\hrule}
\phantom{\frac{\frac{1}{1}}{\frac{1}{1}}} M
\ \ &
\tilde G
{}^{\phantom{\frac{\frac{1}{1}}{\frac{1}{1}}}}
_{\phantom{\frac{\frac{1}{1}}{\frac{1}{1}}}}
&
\rho: \tilde G \to \hat G
{}^{\phantom{\frac{\frac{1}{1}}{\frac{1}{1}}}}
_{\phantom{\frac{\frac{1}{1}}{\frac{1}{1}}}}
\cr \noalign{\hrule}
\text{Elliptic}
&
\smallmatrix
\phantom{a}\\
\left\{\ A \in M_{4,4}(\C)\ :\ A = 
\left(\smallmatrix \smallmatrix a & - b\\ b & a
\endsmallmatrix 
 & 0\\ 
0 & \smallmatrix 1 & 0\\ 0 & 1
\endsmallmatrix \endsmallmatrix \right) \cdot E \ ,\right.\\
\phantom{a}\\
\text{where}\ E\ \text{is}\  I_{4\times 4}\ \text{or}\ 
\left(\smallmatrix \smallmatrix 1 & 0\\ 0 & -1
\endsmallmatrix & 0 \\ 0 & 
\smallmatrix 0 & 1\\1 & 0
\endsmallmatrix \endsmallmatrix \right) \ \\
\phantom{a}\\
\left.\text{and}\ a^2 + b^2 \neq 0\  \right\}\ \simeq\ \C^* \times \Z_2\\
\phantom{a}
\endsmallmatrix
&
\smallmatrix
\phantom{a}\\
\rho\left(\smallmatrix \smallmatrix a & - b\\ b & a
\endsmallmatrix 
 & 0\\ 
0 & \smallmatrix 1 & 0\\ 0 & 1
\endsmallmatrix \endsmallmatrix \right) = 
\left(\smallmatrix a & 
b\\ -b & 
a\endsmallmatrix \right)\\
\phantom{a}\\
\rho\left(\smallmatrix \smallmatrix 1 & 0\\ 0 & -1
\endsmallmatrix & 0 \\ 0 & 
\smallmatrix 0 & 1\\1 & 0
\endsmallmatrix \endsmallmatrix \right)  = \left(\smallmatrix 0 & 1 \\ 
1 & 0\endsmallmatrix \right)
\endsmallmatrix
\cr \noalign{\hrule}
\text{Hyperbolic}
&
\smallmatrix
\phantom{a}\\
\left\{\ A \in M_{4,4}(\R)\ :\ A = 
\left(\smallmatrix \smallmatrix a & 0\\  0 & a
\endsmallmatrix
& 0\\ 
0 & \smallmatrix b & 0\\  0 & b
\endsmallmatrix\endsmallmatrix \right) \cdot E \ ,\right.\\
\phantom{a}\\
\text{where} a, b > 0\ \text{and}\ E\ \text{is}\  I_{4\times 4}\ ,\ 
\left(\smallmatrix \smallmatrix 0 & 1\\  1 & 0
\endsmallmatrix & 0 \\ 
0 & \smallmatrix 1 & 0\\  0 & 1
\endsmallmatrix\endsmallmatrix \right)\ ,\\
\phantom{a}\\
\left.
\left(\smallmatrix \smallmatrix 1 & 0\\  0 & 1
\endsmallmatrix & 0 \\ 
0 & \smallmatrix 0 & 1\\  1 & 0
\endsmallmatrix\endsmallmatrix \right)\ 
\text{or}\  
\left(\smallmatrix 0 & \smallmatrix 1 & 0\\  0 & 1
\endsmallmatrix \\ 
\smallmatrix 1 & 0\\  0 & 1
\endsmallmatrix & 0\endsmallmatrix \right)\  \right\}\ 
\simeq\ \R_* \times \R_* \times \Z_2\\
\phantom{a}
\endsmallmatrix
&
\smallmatrix
\phantom{a}\\
\rho\left(\smallmatrix\smallmatrix a & 0\\ 0 & a\endsmallmatrix & 0\\ 
0 & \smallmatrix b & 0\\ 0 & b
\endsmallmatrix \endsmallmatrix \right) =  \left(\smallmatrix a^2  & 0\\ 
0 & b^2\endsmallmatrix \right)\ ,\ 
\rho\left(\smallmatrix \smallmatrix 0 & 1\\  1 & 0
\endsmallmatrix & 0 \\ 
0 & \smallmatrix 1 & 0\\  0 & 1
\endsmallmatrix\endsmallmatrix \right) =  \left(\smallmatrix -1  & 0\\ 
0 & 1\endsmallmatrix \right)\ \\ 
\phantom{a}\\
\rho\left(\smallmatrix \smallmatrix 1 & 0\\  0 & 1
\endsmallmatrix & 0 \\ 
0 & \smallmatrix 0 & 1\\  1 & 0
\endsmallmatrix\endsmallmatrix \right) = \left(\smallmatrix 1  & 0\\ 
0 & -1\endsmallmatrix \right)\ ,\ 
\rho\left(\smallmatrix 0 & \smallmatrix 1 & 0\\  0 & 1
\endsmallmatrix \\ 
\smallmatrix 1 & 0\\  0 & 1
\endsmallmatrix & 0\endsmallmatrix \right) = \left(\smallmatrix 0  & 1\\ 
1 & 0\endsmallmatrix \right)\\
\phantom{a}
\endsmallmatrix
\cr \noalign{\hrule}
}}
\smallskip
\centerline{\bf Table 2}
\endroster
\endproclaim
\demo{Proof} (1) Let $[a_o] = \left[\matrix 1 & 0\\ 0 & 1 \endmatrix\right]$. By definitions, 
the conormal frames $\theta = (\theta^1, \theta^2) \in \Psi^{-1}([a_o]) \cap  E_x$ 
are exactly  those whose corresponding 
2-forms $(\widetilde{d\theta}^1, 
\widetilde{d\theta}^2)$ represent two space-like vectors of 
an orthonormal basis for the Lorentz metric ${\underset{(\tau)}\to G}$,
 for some fixed choice of the 4-vector $\tau$. Choose
a basis 
 $(\theta^1, \theta^2, e^3, e^4, e^5, e^6)$ of $T^*_xM$, with $e^4|_{\D_x} = J^* e^3|_{\D_x}$ and $e^6|_{\D_x}
 = J^* e^5|_{\D_x}$ and 
with $(e^3\wedge e^4 \wedge e^5 \wedge e^6)(\tau) = 1$, and consider the 
orthonormal basis $\Cal B = (\xi^0, \dots, \xi^3)$ 
of $\Lambda^2_\C\D_x$ given by
$$\xi^0 = \left.\frac{1}{2}\left(e^3 \wedge e^4 + e^5 \wedge e^6\right)\right|_{\D_x}\ ,\qquad 
\xi^1 = \left.\frac{1}{2}\left(e^3 \wedge e^4 - e^5 \wedge e^6\right)\right|_{\D_x}\ ,\tag 3.2$$
$$\xi^2 = \left.\frac{1}{2}\left(e^4 \wedge e^6 + e^3 \wedge e^5\right)\right|_{\D_x}\ ,\qquad 
\xi^3 = \left.\frac{1}{2}\left(e^3 \wedge e^6 - e^4 \wedge e^5\right)\right|_{\D_x}\ .\tag 3.3$$
There always exists an orthonormal basis $\Cal B'$ for $\Lambda^2_\C\D_x$, 
in which the last two vectors are exactly 
the 2-forms $(\frac{1}{2}\widetilde{d\theta}^1, 
\frac{1}{2}\widetilde{d\theta}^2)$. We may also assume that $\Cal B'$ 
has the same orientation and the same time-direction of $\B$. Since, by Lemma 2.3, 
the group $\GL(\D_x, \C)$ acts on $\Lambda^2_\C\D_x$ as the connected component
of the identity of $\operatorname{CO}_{3,1}(\R)$, there is some element 
 $A \in \GL(\D_x, \C)$ which maps $\Cal B'$ into $\Cal B$ and such that $
\xi^2 = \frac{1}{2}A^*\widetilde{d\theta}^1$ and 
$\xi^4 = \frac{1}{2} A^*\widetilde{d\theta}^2$.
This implies that, if we consider any $\Cal A \in GL(T_xM,\R)$ which preserves $\D_x$ and induces 
on $\D_x$ the transformation $A$, then 
the 2-forms 
$\widetilde{d\theta}^i$ can be written in terms of  $e'{}^i = \Cal A^* e^i$ as  in Table 1.\par
Similar  arguments prove  (1)  when 
 $[a_o] = \left[\matrix 0 & 1\\ 1 & 0 \endmatrix\right]$.\par
\medskip
(2) From the $\GL_2(\R)$-equivariance of $\Psi$, it follows immediately that 
$\hat E$ is a reduction of $E$ with structure group equal to 
the stability subgroup $G_{[a_o]} \subset GL_2(\R)$ on $\Bbb P(S_{2\times 2}(\R))$. 
 This stability subgroup is  given in Table 1.\par
\medskip
(3) It can be checked  using just definitions.\qed
\enddemo
\medskip
The relevance of the bundle $\hat E$ comes from the following fact.  
For any (local) diffeomorphism $\phi: M \to M$, let us denote by 
$\hat \phi$ the associated lifted map on 
$$\hat \phi: T^* M \times T^* M \to T^*M \times T^*M \ ,\qquad 
\hat \phi(\theta^1,  \theta^2) \= ((\phi^{-1})^* \theta^1, (\phi^{-1})^*\theta^2)\ .$$
Notice that if $\phi: M \to M$ preserves the distribution $\D$ (i.e.  $\phi_*(\D) \subset \D$)
then the  lifted map $\hat \phi$ is a local diffeomorphism
of the conormal frame bundle $E$ into itself  such that
$\hat \phi^*(\varpi) = \varpi$. 
If we consider a (local) CR diffeomorphism $\phi: M \to M$ (that is  such 
that $\phi_*(\D) \subset \D$ and $\phi^*J = J$), then we have  the following crucial 
property of $\hat E$.\par
\medskip
\proclaim{Proposition 3.2} Let $\hat E \subset E$
be the subbundle defined in Proposition 3.1. 
Then:
\roster
\item"(i)" for any (local) CR diffeomorphism  $\phi: M\to M$, 
the lifted map $\hat \phi$ satisfies  $\hat \phi(\hat E) \subset \hat E$;
\item"(ii)" a  (local) diffeomorphism $\varphi: \hat E \to \hat E$ 
 is  the lifted map $\varphi = \hat \phi$
of some  (local) 
diffeomorphism $\phi: M \to M$ preserving $\hat E$ (not necessarily CR),
if and only if it  satisfies 
$$\varphi^*(\varpi|_{\hat E}) = \varpi|_{\hat E}\ .\tag 3.4$$
\endroster
\endproclaim
\demo{Proof} Let $(\theta^1, \theta^2) \in \hat E_x$, for some $x\in M$, and let $(\theta^1, \theta^2, 
e^3, \dots, e^6)$ be a basis of $T^*_xM$ as in Proposition 3.1 (1). If $\phi: M\to M$ 
is a CR diffeomorphism defined on a neighborhood of $x$, one can check 
that the 1-forms
$\theta^1{}' = \phi^*(\theta^1)$,  $\theta^2{}' = \phi^*(\theta^2)$, $e^i{}' = \phi^*(e^i)$, 
$i = 3, \dots, 6$, constitute a basis for $T^*_{\phi(x)} M$ and the 2-forms
$\widetilde{d\theta^1{}'}$ and $\widetilde{d\theta^2{}'}$ are written in terms of 
$e^i{}'$ as in Proposition 3.1 (1). From this, it follows  that $\Psi(\theta^1{}',
\theta^2{}') = \Psi(\theta^1,
\theta^2) = [a_o]$ and  that $\hat \phi(\theta^1, \theta^2) \in \hat E$. The 
claim that $\hat \phi(\varpi|_{\hat E}) = \varpi|_{\hat E}$
follows directly from definitions.\par
The proof that a  (local) diffeomorphism $\varphi: \hat E \to \hat E$  is the lifted map 
of some CR transformation of $M$  if and only if (3.4) holds 
can be obtained by the same line of arguments of Proposition VI.1.3 in \cite{13}. \qed 
\enddemo
\remark{Remark 3.3} Notice that, in addition to $(\D,J^{(0)}) = (\D,J)$ and 
$(\D, J^{(1)}) = (\D, - J)$,  the manifold $M$ carries  
two  natural (almost) 
CR structures $(\D, J^{(i)})$, $i = 2,3$, defined as follows. \par 
At a point $x\in M$, pick any  pair $\theta = (\theta^1, \theta^2) \in \hat E_x$ 
and some corresponding 1-forms $(e^3, e^4, e^5, e^6)$ which satisfy  Proposition
3.2 (1). Then $\D_x$ splits into the $J$-invariant 
subspaces  $\D^1_x = \ker e^5|_{\D_x} \cap 
\ker e^6|_{\D_x}$ and $\D^2_x = \ker e^3|_{\D_x} \cap 
\ker e^4|_{\D_x}$. This splitting is independent of 
the choice of $\theta$ and of the $e^i$'s (see \cite{5}). 
The CR structures  $(\D, J^{(2)})$ and $(\D,J^{(3)})$
 are defined by 
$$J^{(2)}_x = J|_{\D^1_x} \oplus (- J|_{\D^2_x})\ ,\qquad J^{(3)}_x = 
(-J|_{\D^1_x}) \oplus (J|_{\D^2_x})\ ,\qquad x\in M\ .$$ 
It is not difficult to check that  any
local transformation
$\phi: M \to M$, 
which preserves $\D$ and such that $\phi_*(J^{(i)}) = J^{(j)}$ for 
some $0\leq i,j \leq 3$, 
has a lift which maps $\hat E$ into $\hat E$.
\endremark
\bigskip
\bigskip
\subhead 4. The Chern-Moser bundle of an elliptic or hyperbolic manifold
\endsubhead
\bigskip
In the following, we will continue adopting the notation of the previous section. 
In particular, $\varpi = (\varpi^1, \varpi^2)$ will be used also to denote 
the restriction of the tautological pair of  $E = E(M,\D)$ 
on the reduction $\hat E \subset E$. \par
Moreover, since $E$ is a principal bundle over $M$, with structure 
group $\GL_2(\R)$, we may  consider the fundamental vector fields
$E^i_j{}^*$ on $E$, which are 
determined by the elements $E^i_j = (\delta^i_j) \in \goth{gl}_2(\R)$. 
Recall the $E^i_j{}^*$ are vertical vector fields which span at all points
the vertical distribution of  $E(M,\D)$. 
Finally, we adopt also the following notation.\par
\medskip
\remark{Notation 4.1} We will use  latin letters $a, b, c, d$ to denote  indices
which run between $1$ and $2$; we will use the letters $i, j, k, \ell$
to denote indices which run between $3$ and $6$; we will use greek
letters $\alpha, \beta, \gamma, \delta$ to denote indices which run between $7$ and $8$. \par
With capital latin letters $I, J, K, L$, we will denote indices which may run between 
$1$ through $8$.\par
Consider a frame $(e_1, \dots, e_8) \subset T_\theta\hat E$ 
at a point of $\hat E$ and 
denote  by $(e^1, \dots, e^8) \subset T^*_\theta\hat E$
the associated dual frame. Then we will  use the symbols
$e^1_1, e^2_1, e^1_2, e^2_2$ to denote a quadruple of 1-forms, 
which depends
on the forms $e^7, e^8$ according to the following rules: 
\roster
\item"i)" if $M$ is elliptic, 
we set  $e^1_1 = e^2_2 = e^7$ and $- e^1_2 = e^2_1 = e^8$; 
\item"ii)" if $M$ is hyperbolic, we set
$e^1_1 = e^7$, $e^2_2 = e^8$ and $e^1_2 = e^2_1 = 0$.
\endroster
Furthermore,  we  define the following 
$\C$-valued 1-forms:
$E^0 = e^1 + i e^2$, $E^1 = e^3 + i e^4$, $E^2 = e^5 + i e^6$, $E^0_0 = e^1_1 + i e^2_2$, 
$\Omega = \varpi^1 + i \varpi^2$.
\endremark
\bigskip
We have now all ingredients to introduce the concept of ``adapted frames" and of ``Chern-Moser
bundle of an elliptic or hyperbolic manifold".\par\
\medskip
\definition{Definition 4.2}
A frame $(e_1, \dots, e_N) \subset T_\theta \hat E$ 
 is  called {\it adapted to the CR structure\/} if the following conditions
are satisfied:
\roster
\item"1)" the vectors  $e_7$ and $e_8$ are  equal to the  vectors $\tilde e_7$, $\tilde e_8$
defined by
$$\tilde e_7 = \left\{ \matrix E^1_1{}^*|_\theta + E^2_2{}^*|_\theta & \text{if}\ M\ \text{is elliptic}\ ,\\
E^1_1{}^*|_\theta  & \text{if}\ M\ \text{is hyperbolic}\ ,
\endmatrix\right.\tag{$4.1_1$} $$
$$\tilde e_8 = \left\{ \matrix - E^1_2{}^*|_\theta + E^2_1{}^*|_\theta & \text{if}\ M\ \text{is elliptic}\ ,\\
E^2_2{}^*|_\theta  & \text{if}\ M\ \text{is hyperbolic}\ ,
\endmatrix\right.\tag{$4.1_2$} $$
\item"2)" the vectors $e_a$, $a = 1,2$,  are such that $\varpi^a(e_b) = \delta^a_b$; 
\item"3)"  the vectors $e_i$, $i = 3,4,5,6$,  satisfy 
$J \hat \pi_*(e_3) = \pm \hat \pi_*(e_4)$,  $J \hat \pi_*(e_5) = \pm \hat \pi_*(e_6)$
and 
the linear equations
$$\varpi^a(e_i) = 0 \ ,\qquad 
\ d\varpi^a(e_i, e_\alpha) = 0\tag 4.2$$
for any $a = 1,2$ and $\alpha = 7,8$ plus the following conditions
 (here $\Omega = \varpi^1 + i \varpi^2$, 
$E_0 = \frac{1}{2}(e_1 - i e_2)$,
$E_1 = \frac{1}{2}(e_3 - i e_4)$ and $E_2 = \frac{1}{2}(e_5 - i e_6)$):
\medskip
\moveright 0.5 cm
\vbox{\offinterlineskip
\halign {\strut\vrule\hfil\ $#$\ \hfil
 &\vrule\ \vrule\hfil\ $#$\ \hfil
\vrule\cr
\noalign{\hrule}
\text{Elliptic}
& 
{\left\{\matrix  
{d\Omega(E_1, E_0) = d\Omega(\bar E_1, \bar E_0) =
d\Omega(E_2, \bar E_0) = d\Omega(\bar E_2,  E_0) = 0
}^{\phantom{\frac{\frac{1}{1}}{\frac{1}{1}}}}
\\
{d\Omega(\bar E_1, E_0) = d\Omega(E_2, E_0) = 0 
}^{\phantom{\frac{\frac{1}{1}}{\frac{1}{1}}}}
\\
{d\Omega(\bar E_1, E_2)  = 1\qquad d\Omega(E_1, \bar E_2) = 0}^{\phantom{\frac{\frac{1}{1}}{\frac{1}{1}}}}\\
{d\Omega(E_1, \bar E_1) = d\Omega( E_2, \bar E_2) = 
d\Omega(E_1, E_2) = d\Omega(\bar E_1, \bar E_2) = 0}^{\phantom{\frac{\frac{1}{1}}{\frac{1}{1}}}}
_{\phantom{\frac{\frac{1}{1}}{\frac{1}{1}}}}
\endmatrix
\right.
}^{\phantom{\frac{\frac{1}{1}}{\frac{1}{1}}}}
_{\phantom{\frac{\frac{1}{1}}{\frac{1}{1}}}}
\cr \noalign{\hrule}
\text{Hyperbolic}
&
{\left\{\matrix  
{d\varpi^a(e_i, e_a) =0\ \text{for any}\ a = 1,2}^{\phantom{\frac{\frac{1}{1}}{\frac{1}{1}}}}\\
{d\varpi^1(e_i, e_2) =0 \ \ \text{if}\ i = 3,4\ \text{and}\ d\varpi^2(e_i, e_1) =0
\ \ \text{if}\ i = 5,6}^{\phantom{\frac{\frac{1}{1}}{\frac{1}{1}}}}
\\
d\varpi^1(e_i, e_j)  = {\left\{\matrix 1  & \text{if}\  i = 3, j = 4 \\
0 & \text{if}\ \{i,j\} \neq \{3,4\}\endmatrix \right.}^{\phantom{\frac{\frac{1}{1}}{\frac{1}{1}}}} \\
d\varpi^2(e_i, e_j)  = {\left\{\matrix 1  & \text{if}\ i = 5, j = 6 \\
0  & \text{if}\ \{i,j\} \neq \{5,6\}\endmatrix \right.}^{\phantom{\frac{\frac{1}{1}}{\frac{1}{1}}}}
_{\phantom{\frac{\frac{1}{1}}{\frac{1}{1}}}}
\endmatrix
\right.
}^{\phantom{\frac{\frac{1}{1}}{\frac{1}{1}}}}
_{\phantom{\frac{\frac{1}{1}}{\frac{1}{1}}}}
\cr \noalign{\hrule}
}}
\smallskip
\centerline{\bf Table 3}
\endroster
\enddefinition
\bigskip
The conditions given in  Definition 4.2 can be totally reformulated into 
conditions on the dual coframe of an adapted frame.  Such conditions 
are the following.\par
\medskip
\proclaim{Lemma 4.3}ÊA frame $(e_1, \dots, e_8)$ of a tangent space $T_\theta \hat E$ 
is a adapted to the CR structure if and only if the dual coframe
$(e^1, \dots, e^8) \subset T^*_\theta \hat E$
 satisfies the following 
conditions: 
\roster
\item"i)" for $\alpha = 7,8$, the 1-forms $e^\alpha$  satisfy 
$e^\alpha(\tilde e_\beta) = \delta^\alpha_\beta$,
where the $\tilde e_\beta$, with $\beta = 7,8$, are the vectors defined in (4.1);
\item"ii)" for $a = 1,2$, $e^a = \varpi^a|_\theta$; 
\item"iii)" for $i = 3,4,5,6$, the 1-forms $e^i$ vanish on any of the vectors 
$\tilde e_\beta$, $\beta = 7,8$,  given in  (4.1), and there exists
a coframe $(\theta^1, \dots, \theta^6)$ in $T^*_{x} M$, $x = \hat \pi(\theta)$, such that 
$\hat \pi^*(\theta^i) = e^i$, $i = 3,\dots, 6$, and 
$J^*\theta^3|_{\D_x} = \pm \theta^4|_{\D_x}$, $J^*\theta^5|_{\D_x} = \pm \theta^6|_{\D_x}$; 
\item"iv)"  the differentials $d\varpi^a$, $a = 1,2$, evaluated at the point $\theta$, are equal 
to the following expressions for some suitable 
constants $s,s',t,t', \sigma, \tau$
(here $\Omega = \varpi^1 + i \varpi^2$, 
$E^1 = e^3 + i e^4$, $E^2 = e^5 + i e^6$, $E^0_0 = e^7 + i e^8 = e^1_1 + i e^2_2$): 
\medskip
\moveright 1 cm
\vbox{\offinterlineskip
\halign {\strut\vrule\hfil\ $#$\ \hfil
 &\vrule\ \vrule\hfil\ $#$\ \hfil
\vrule\cr
\noalign{\hrule}
\phantom{\frac{\frac{1}{1}}{\frac{1}{1}}}\D
\ \ &
d\varpi^a
\cr \noalign{\hrule}
\text{Elliptic}
& d\Omega  +  E^0_0 \wedge \Omega = 
\bar E^1 \wedge E^2 + \sigma\ E^1 \wedge \bar \Omega + 
\tau\ \bar E^2 \wedge \bar \Omega
{}^{\phantom{\frac{\frac{1}{1}}{\frac{1}{1}}}}
_{\phantom{\frac{\frac{1}{1}}{\frac{1}{1}}}}
\cr \noalign{\hrule}
\text{Hyperbolic}
&
\matrix
\ \\ 
d\varpi^1 + 
e^1_a \wedge \varpi^a = e^3 \wedge e^4 + 
s\ e^5 \wedge \varpi^2 + t\ e^6 \wedge \varpi^2\\
\ \\
d\varpi^2 + 
e^2_a \wedge \varpi^a = e^5 \wedge e^6 + 
s'\ e^3 \wedge \varpi^1 + t'\ e^4 \wedge \varpi^1\\
\ 
\endmatrix 
{}^{\phantom{\frac{\frac{1}{1}}{\frac{1}{1}}}}
_{\phantom{\frac{\frac{1}{1}}{\frac{1}{1}}}}
\cr \noalign{\hrule}
}}
\smallskip
\centerline{\bf Table 4}
\endroster
\endproclaim
\demo{Proof} Let us  check that if 
a coframe $(e^1, \dots, e^8)$ satisfies  i) - iv), then the dual basis 
$(e_1, \dots, e_N)$ is an 
 adapted frame.\par
First of all, from i), ii) and iii), it is clear that the vectors 
$e_\alpha$, with $\alpha = 7,8$, must coincide with the vectors $\tilde e_\alpha$
of (4.1). Also   (2) of Definition 4.2 is immediately satisfied 
by the vectors $e_1$ and $e_2$.  Finally, if we consider the remaining vectors 
$e_i$, $i = 3,4,5,6$, and we plug them into the tautological 1-forms $\varpi^a|_\theta = e^a$
and into  $d\varpi^a|_\theta$, we see that all conditions of 
Definition 4.2 (3) are  satisfied.\par
Conversely, assume that $(e^1, \dots, e^8)$ is a coframe, which is dual to an 
adapted frame. It is clear that i) and iii) are satisfied. Also ii) is satisfied, since $\varpi^a|_\theta$, 
$a = 1,2$, gives the value $1$ if and only if it is evaluated to the vector $e_a$. 
Finally, from  (3) of Definition 4.2, it follows that, modulo terms of type
 $\varpi^a\wedge e^\alpha$, $a =1,2$, $\alpha = 7,8$, the expressions for 
$d\varpi^a|_\theta$ have to be as in Table 4.  Then, using the action of the vector fields $E^i_j{}^*$ on the 
tautological 1-forms $\varpi^a$, one can compute 
the values $d\varpi^a((E^i_j)^*,e_a)$ and  check directly that 
the terms of type $\varpi^a\wedge e^\alpha$, $a =1,2$, $\alpha = 7,8$, 
appearing in the expressions for
$d\varpi^a|_\theta$,  are those given  in  Table 4 (at this regard, 
see also the proof of next Lemma 4.4).\qed
\enddemo
\medskip
It is fair to ask if there exists at least one adapted frame at any  point of $\hat E$. 
The answer is yes as it is proved in the following lemma.
\medskip
\proclaim{Lemma 4.4} There exists at least 
one adapted frame at any  $\theta_o = (\theta^1_o, \theta^2_o) \in \hat E$. 
\endproclaim
\demo{Proof} Let  $x$ be the point $x = \hat \pi(\theta_o) \in M$ and 
consider a local section $\hat \theta: \Cal U \subset M \to \hat E$ such 
that $\hat \theta_x = \theta_o$. Any element $\theta \in \hat \pi^{-1}(\Cal U) 
\subset E(M, \D)$  can be written as
$$\theta = A^{-1}\cdot \hat \theta_y$$
where $y = \pi(\theta) \in \Cal U$ and $A \in \GL_2(\R)$. Moreover, if we consider the 1-form 
$\hat \varpi = \hat \pi^* \hat \theta$, we have that the tautological pair
$\varpi$ can be written at a point $\theta = A^{-1}\cdot \hat \theta_y$ as
$$\varpi|_{A^{-1}\cdot \hat \theta_y} = A^{-1} \cdot \hat \varpi|_{A^{-1}\cdot \hat \theta_y}\ .$$
Therefore, by the fact that $\theta_o = \hat \theta_x$, we get that 
 $d\varpi$ at $\theta_o$ can be written as
$$d\varpi|_{\theta_o} = - (dA) \wedge \varpi|_{\theta_o} + d\hat \varpi|_{\theta_o}\ .\tag 4.3$$
It is not difficult to realize that $(dA)$ is a $2\times 2$-matrix of 1-forms
$$(dA) =  \left(\matrix e^1_1 & e^1_2\\
e^2_1 & e^2_2 \endmatrix \right)$$
with entries $e^i_j$ which satisfy the linear relations
described in i) or ii) in 
Notation 4.1.\par
On the other hand, $d\hat \varpi|_{\theta_o} = \hat \pi^* d\hat \theta_x$ and, by 
Proposition 3.1, we can find  1-forms $\hat e^3, \dots, \hat e^6$ in 
$T^*_xM$ so that $d\hat \theta^a_x$ can be written as
$$d\hat \theta^a_x = H^a_{ij} \hat e^i\wedge \hat e^j + M^a_{i b} \hat e^i\wedge \theta^b_o
+ N^a_{b c} \theta^b_o \wedge \theta^c_o\ ,\tag 4.4$$
where the constants $H^a_{ij}$ are so that $\widetilde{d\theta_o} = 
d\hat \theta|_{\D_x \times \D_x}$ is in one of the forms listed in Table 1.
So, if we set 
$$e^a = \hat \pi^*\theta^a_o\ , \qquad e^i = \hat \pi^* \hat e^i\ ,$$
we get that  the 
 coframe  $(e^1, \dots, e^6, e^i_j)$ satisfies i), ii) and iii) of Lemma 4.3. Moreover, 
the differentials $d\varpi^a|_{\theta_o}$  are of the form
$$d\varpi^a|_{\theta_o} + e^a_b \wedge \varpi^b = 
H^a_{ij} e^i\wedge e^j + M^a_{i b} e^i\wedge \varpi^b
+ N^a_{b c} \varpi^b \wedge \varpi^c\ ,\tag 4.5$$
for some real numbers $M^a_{ib}$ and $N^a_{ib}$. 
Now, replacing the 1-forms $e^i, e^a_b$ with other 1-forms of the kind $e^i + {\Cal A}^i_a e^a$
and $e^a_b + {\Cal B}^a_{bc} e^c + {\Cal C}^a_{bi} e^i$ for some constants  
${\Cal A}^i_a$, ${\Cal B}^a_{bc}$, ${\Cal C}^a_{bi}$, one gets another  coframe which still 
satisfies i), ii) and iii) of Lemma 4.3. Moreover, 
if the coefficients ${\Cal A}^i_a$ and 
${\Cal B}^i_{ja}$ are suitably chosen, one can obtain that 
several of the constants $ M^a_{i b}$ and $N^a_{b c}$, which appear in relation  (4.5) for this 
new coframe, are equal to $0$. Choosing the coefficient, so 
that a  maximal  number of the constants 
 $ M^a_{i b}$ and $N^a_{b c}$ vanishes, one obtains the equations of Table 4.
\qed
\enddemo
\bigskip
\definition{Definition 4.5} Let 
$M$ be an elliptic  or 
hyperbolic manifold. 
The {\it extended Chern-Moser bundle 
of  $M$\/} is the set $\hat \P$  of all adapted frames of $\hat E$. 
The {\it  
Chern-Moser bundle of $M$\/} is the subset $\P \subset \hat \P$ given by all adapted frames
such that $J \hat \pi_*(e_3) = + \hat \pi_*(e_4)$ and $J \hat \pi_*(e_5) = + \hat \pi_*(e_6)$.\par
The set $\hat \P$ has a natural 
structure of fiber bundle over $M$ given by the projection 
$$\pi = \hat \pi \circ \pi_o: \P\to M\ ,\qquad \hat \pi: \hat E \to M\ ,\ \pi_o: \P \to \hat E\ .$$
The {\it tautological 1-form of $\P$\/}\ is the 6-tuple 
 $\omega = (\omega^1, \dots, \omega^6)$, where the $\omega^i$'s are the 1-forms 
defined by 
$$\pi_*(X) = \sum_{i = 1}^6 \omega^i(X)\cdot \hat \pi_*(e_i)\ \text{or, equivalently,}\ 
\omega^i(X) = e^i(\pi_o{}_*(X))\ ,\ i = 1, \dots, 6\ , $$
for any $X\in T_u\hat \P$ at a frame $u = (e_1, \dots, e_8) \in \hat\P$.  
\enddefinition
\bigskip
\bigskip
\subhead 5. The Chern-Moser bundle $\P$ is a principal bundle over $M$
\endsubhead
\bigskip
The aim of this section is the proof of Theorems 5.2 and 5.3 below, which claim that
both $\hat \P$ and 
$\P$ admit a natural structure of  principal bundle over $M$.
Their proofs 
require a preliminary result, given in the next 
Proposition 5.1.\par 
\medskip
\subsubhead 5.1 The natural action of the structure group of $\hat \pi: \hat E \to M$  on
$\hat \P$  
\endsubsubhead\par
\medskip
The first step for Theorem 5.3 consists in showing 
that the structure group $\hat G$ of $\hat \pi: \hat E \to M$ admits 
a lifted action on $\pi_o: \hat \P \to \hat E$. Namely, \par
\medskip
\proclaim{Proposition 5.1} Let  $\hat G$ be the
structure group of $\hat \pi: \hat E \to M$ as described in Proposition 3.1. 
Then there exists a  right
action of $\hat G$ on $\hat \P$, which commutes with any diffeomorphism of 
$\hat \P$ induced by a CR transformation of $M$.
\endproclaim
\demo{Proof} 
In the following, for any element $A \in \hat G$,  we will denote by $\tilde A = \rho^{-1}(A)\in \tilde G$ the 
corresponding element in the group $\tilde G \subset \GL_2(\C)$ given in Table 2. Moreover, 
 for any $A\in \hat G$, let $R_A: \hat E \to \hat E$ be the right action of $A$, i.e. 
$$R_A(\theta^1, \theta^2) = ((A^{-1})^1_a \theta^a, (A^{-1})^2_b \theta^b) \=  (A^{-1}) \cdot \theta$$
Finally, for any frame $u = (e_1, \dots, e_8) \subset T_\theta \hat E$, 
let $\Psi_A(u) = (e'_1, \dots, e'_8)$ be the frame 
$$e'_a =  R_A{}_*(A_a^b e_b)\ ,\quad e'_i = 
R_A{}_*(\tilde  A_i^j e_j)\ ,
\quad e'_\alpha = R_A{}_*(e_\alpha) \tag 5.1$$
(we  adopt the convention on indices 
of Notation 4.1).
It is simple to check that, for any two elements $A, A' \in \hat G$
$$\Psi_A \circ \Psi_{A'} = \Psi_{A' \cdot A}$$
and  hence that the map $A \mapsto \Psi_A$ gives 
a  right action of $\hat G$ on the space of linear frames
of $\hat E$. If we  show that this action maps adapted frames into adapted frames, we are done. \par
Assume that $(e_1, \dots, e_8)$ is an adapted frame and let 
$(e^1, \dots, e^8)$ be the corresponding dual coframe. Since the fundamental vector fields determined 
by the structure group $\hat G$ are mapped into itself by any diffeomorphism 
$R_A$, with $A \inÊ\hat G$, it is clear that the vectors
$e'_\alpha$, $\alpha = 7,8$,  satisfy (1) of Definition 4.2 whenever the vectors $e_\alpha$'s do. 
Notice also that, for any $X \in T_\theta \hat E$
$$R_A^*(\varpi^a)(X)|_\theta = \varpi^a_{A^{-1}\cdot \theta}(R_A{}*(X)) = $$
$$ = 
((A^{-1})^a_b \theta^b)(\hat \pi_* \circ R_A{}_*(X)) = ((A^{-1})^a_b \theta^b)(\hat \pi_*(X)) = 
((A^{-1})^a_b\varpi^b)(X)|_{\theta}\ .\tag 5.2$$
Therefore, for $a,b, c = 1,2$, 
$$\varpi^a(e'_b) = A_b^c R_A{}^*(\varpi^a)(e_c) = 
A^c_b (A^{-1})^a_d  \delta^d_c  = \delta^a_b\ ,$$
and hence also  (2) of Definition 4.2 is satisfied. 
Similarly,  using (5.2) and the fact that the vectors $e_i$'s satisfy (4.2),
one can check that also the vectors 
$e'_i = \tilde A_i^j R_A{}_*(e_j)$ satisfy (4.2).\par
It remains to check if the vectors $e'_i$ satisfy  the conditions of Table 1. But this can
be done, just using the explicit expressions of 
the matrices $A^a_b\in \hat G$ and $\tilde A^i_j \in \tilde G$ and  the fact that, by construction, 
$d\varpi^a|_{A^{-1}\cdot \theta}(e'_j, e'_b) = 
\tilde A^i_j (A^{-1})^a_c A^d_b d\varpi^c(e_i, e_d)$ 
for any $a,b = 1,2$ and any $i,j = 3, 4,5,6$. \qed
\enddemo
\bigskip
\subsubhead 5.2 The bundle $\pi: \hat\P\to \hat E$ is a principal bundle (elliptic case)
\endsubsubhead\par
\medskip
Assume that  $M$  is elliptic. Recall that, for any adapted frame 
$(e_1, \dots, e_8)$ at $\theta \in \hat E$,  with corresponding dual coframe 
$(e^1, \dots, e^8)$ the following holds (see Notation 4.1)
$$d\Omega + E^0_0 \wedge \Omega = \bar E^1 \wedge E^2 + \sigma E^1 \wedge \bar \Omega
+ 
\tau \bar E^2 \wedge \bar \Omega\ ,\tag 5.3$$
where $\sigma$ and $\tau$ are two suitable complex numbers, depending on the 
frame $(e_i)$. \par
Now, consider a new 
coframe $(e'{}^i)$, which satisfies i), ii) and iii) of Lemma 4.3. Assume also that, 
if $(e^i)$ satisfies (iii) of Lemma 4.3 with 
$J^*\theta^3= \epsilon \theta^4$ and $J^* \theta^5 = \epsilon' \theta^6$, for 
some  $\epsilon, \epsilon' = \pm 1$, 
then the new coframe satisfies (iii) of Lemma 4.3 with the same signs $\epsilon$ and 
$\epsilon'$. Then $(e'{}^i)$ 
is of the form
$$e^a{}' = e^a \quad (\ \Rightarrow\ \Omega' = 
e'{}^1 + i e'{}^2 = e^1 + i e^2 = \Omega\ )\ ,\tag5.4$$
$$E^i{}' = C^i_j E^j + B^i (e^1 + i e^2) + B^i_b (e^1 - i e^2) = 
C^i_j E^j + B^i \Omega + B^i_b \bar \Omega\ ,\tag 5.5$$
$$E^0_0{}'= E^0_0 + \hat B_j E^j + \tilde B_j \bar E^j +  \hat A \Omega + \tilde A \bar \Omega\ ,
\tag 5.6$$
where $C^i_j$, $B^i_a$, $\hat B_j$, $\tilde B_j$, $\hat A$ and $\tilde A$ are 
complex numbers. In order to make $(e'{}^i)$ to satisfy also iv) of Lemma 4.3, we have to require that 
the constants $C^i_j$, $B^i_a$, $\hat B_j$, $\tilde B_j$, $\hat A$ and $\tilde A$  
satisfy some additional conditions. In fact,   plugging  (5.4) - (5.6) into 
 (5.3), we get that 
$$d\Omega = - E^0_0 \wedge \Omega - \hat B_j E^j \wedge \Omega - 
\tilde B_j \bar E^j \wedge \Omega + \tilde A \Omega \wedge \bar \Omega +$$
$$ + \overline{C^1_1} C^2_2\bar E^1 \wedge E^2 + \overline{C^1_2} C^2_2\bar E^2 \wedge E^2 +
\overline{C^1_1} C^2_1\bar E^1 \wedge E^1 + \overline{C^1_2} C^2_1\bar E^2 \wedge E^1 + $$
$$ +   
\overline{C^1_i} B^2 \bar E^i \wedge \Omega + \overline{C^1_i} B^2_b \bar E^i \wedge \bar \Omega
-
\overline{B^1} C^2_i E^i \wedge \bar \Omega - \overline{B^1_b} C^2_i E^i \wedge \Omega +$$
$$ +  
\sigma' C^1_i E^i \wedge \bar \Omega
+ \sigma' B^1 \Omega \wedge \bar \Omega  +  
\tau' \overline{C^2_i} \bar E^i \wedge \bar \Omega  +
\tau' \overline{B^2_b} \Omega\wedge \bar \Omega\ .\tag 5.7$$
By comparison of (5.7) with (5.3), we find  the conditions
$$C^1_1 \overline{C^2_2} = 1\ ,\quad C^1_2 = C^2_1 = 0\ ,\quad
B^1 = B^2_b = 0\ ,
\tag 5.8$$
$$\tilde B_1 = \overline{C^1_1} B^2\ ,\quad 
\hat B_2 = - \overline{B^1_b} C^2_2\ ,\quad \hat B_1 = \tilde B_2 = \tilde A = 0\ .\tag 5.9$$
Moreover, we see that the values  $\sigma, \tau$ associated with the frame $(e_i)$
and the values $\sigma', \tau'$ associated with the frame $(e_i)$ are related by
$$\sigma' = (C^1_1)^{-1} \sigma = \overline{C^2_2} \sigma\ ,\qquad \tau' = (\overline{C^2_2})^{-1} \tau = 
C^1_1 \tau\ .\tag 5.10$$
From such observations, we obtain that $(e'{}^i)$
is a new adapted frame if and only if it is of the form
$$E^0{}' (\= e'{}^1 + i e'{}^2) = E^0\ , \quad 
E^1{}' = \frac{1}{\bar C} E^1 +  \bar F \overline{E^0} \ ,\quad E'{}^2 = C E^2 + H E^0\ ,$$
$$E^0_0{}' = E^0_0 - C F E^2  - \frac{1}{C} H \bar E^1 + A E^0\ .$$
\smallskip
With similar arguments, it can be checked that there is 
no adapted coframe $(e'{}^i)$, which satisfies  (iii) of Lemma 4.3 with signs
$J^*\theta^3= + \epsilon \theta^4$ and $J^* \theta^5 = - \epsilon'\theta^6$ or $J^*\theta^3= - 
\epsilon \theta^4$ and 
$J^* \theta^5 = + \epsilon'\theta^6$.\par
\smallskip
Finally, using  the same arguments of before, we get that any adapted 
coframe $(e'_i)$, which satisfies  iii) of Lemma 4.3 with
$J^*\theta^3= - \epsilon\theta^4$ and $J^* \theta^5 = - \epsilon'\theta^6$, has to be of the 
form 
$$E^0{}' = E^0\ , \quad 
E^1{}' = - \frac{1}{\bar C} \bar E^2 +  \bar F \overline{E^0} \ ,\quad E'{}^2 = C \bar E^1 + H E^0\ ,$$
$$E^0_0{}' = E^0_0 - \frac{1}{C} H E^2 - C F \bar E^1 + A E^0\ .$$
It follows that $\hat \P$ is a principal bundle 
over $\hat E$, with structure group $G$, which is isomorphic to the following group 
of matrices, associated with the transformations of
 the quadruple $(E^0, \overline{E^1}, E^2, E^0_0)$
into the quadruple $(E'{}^0, \overline{E'{}^1}, E'{}^2, E'{}^0_0)$:
$$G = G^{(1)} \cdot G^{(2)} \qquad G^{(1)} = \left\{ \ \left( \matrix 
1 & 0 & 0  & 0 \\
F &  \frac{1}{C} & 0 & 0  \\
H  & 0 & C & 0 \\
A & \frac{H}{C} & - CF & 1
 \endmatrix \right)\ 
A,C,F,H\in \C\ ,\ \right\}\ ,$$
$$
G^{(2)} = \left\{ \ I_{4\times 4}\ ,\ \left( \matrix 
1 & 0 & 0  & 0 \\
0 &  0 & 1 & 0  \\
0  & -1 & 0 & 0 \\
0 & 0 & 0 & 1
 \endmatrix \right)\ \ ,\ \left( \matrix 
1 & 0 & 0  & 0 \\
0 &  - 1 & 0 & 0  \\
0  & 0 & -1 & 0 \\
0 & 0 & 0 & 1
 \endmatrix \right)\ ,\ \left( \matrix 
1 & 0 & 0  & 0 \\
0 &  0 & -1 & 0  \\
0  & 1 & 0 & 0 \\
0 & 0 & 0 & 1
 \endmatrix \right)\ \right\}
\tag 5.11$$
In particular, {\it $\hat \P$ is a smooth manifold\/}. 
Note that $G^{(1)}$ is isomorphic to a subgroup of the isotropy $\HQ$, where $\Q$ is the elliptic 
quadric. \par
\bigskip
\subsubhead 5.3 The bundle $\pi: \hat \P\to \hat E$ is a principal bundle (hyperbolic case)
\endsubsubhead
\par
\medskip
Assume  now  that   $M$  is
hyperbolic. 
For any adapted frame $(e_1, \dots, e_8)$ at $\theta \in \hat E$,  with corresponding dual coframe 
$(e^1, \dots, e^8)$, from Table 4 and adopting Notation 4.1, we may write that 
$$d\varpi^1|_\theta = \frac{i}{2} E^1 \wedge \bar E^1 + \sigma E^2 \wedge e^2 + 
\bar \sigma \bar E^2\wedge e^2 - e^1_1 \wedge e^1\ ,\tag 5.12$$
$$d\varpi^2|_\theta = \frac{i}{2} E^2 \wedge \bar E^2 + \tau E^1 \wedge e^1 + 
\bar \tau \bar E^1\wedge e^1 - e^2_2 \wedge e^2\ ,\tag 5.13$$
for some $\sigma, \tau \in \C$. 
Assume also that
$(e^i)$ satisfies (iii) of Lemma 4.3 with 
$J^*\theta^3= \epsilon \theta^4$ and $J^* \theta^5 = \epsilon' \theta^6$, for 
some fixed values $\epsilon, \epsilon' = \pm 1$. Now, from definitions,  
a new coframe $(e'{}^i)$,
which satisfy (i), (ii) and (iii) of Lemma 4.3, with the same signs $\epsilon$, 
$\epsilon'$ in the equations $J^*\theta'{}^3= \epsilon \theta'{}^4$ and $J^* \theta'{}^5 = 
\epsilon' \theta'{}^6$, 
can be obtained from $(e^i)$
by means of a linear transformation of the following form:
$$e^a{}' = e^a\ ,\qquad E^i{}' = C^i_j E^j + B^i_a e^a\ ,\tag 5.14$$
$$e^1_1{}'= e^1_1 + \hat B^1_j E^j + \overline{\hat B^1_j} \bar E^j +  A^1_a e^a\ ,
\quad
e^2_2{}'= e^2_2 + \hat B^2_j E^j + \overline{\hat B^2_j} \bar E^j +  A^2_a e^a\ ,\tag 5.15$$
for some  $C^i_j, B^i_a, \hat B^i_j \in \C$ and $A^i_a \in \R$. Assuming that  the coframe 
$(e'{}^i)$ satisfies also (iv) (and hence (5.12) and (5.13)),  the following
conditions have to be satisfied
$$|C^1_1|^2 = |C^2_2| = 1\ ,\quad C^a_j = 0 \ \text{if}\ a \neq j\ ,$$
$$B^1_2 = \hat B^1_2 = A^1_2 = 0\ ,\qquad 
B^2_1 = \hat B^2_1 =  A^2_1 = 0\ ,$$
$$\hat B^1_1 = \frac{i}{2} \bar B^1_1 C^1_1\ ,\ \hat B^2_2 = \frac{i}{2} \bar B^2_2 C^2_2\ ;$$
there is no restriction on the coefficients $B^i_i \in \C$ and $A^1_1 \in \R$. \par
\smallskip
On the other hand, 
a similar line of arguments shows that if 
a new coframe $(e'{}^i)$ satisfies (i), (ii) and (iii) of Lemma 4.3, with the  additional conditions
 $J^*\theta'{}^3= - \epsilon \theta'{}^4$ or $J^* \theta'{}^5 = 
-\epsilon' \theta'{}^6$, then the condition (iv) of Lemma 4.3 can never be satisfied and 
hence $(e'{}^i)$ cannot be an adapted coframe.\par
\smallskip
Since this result is independent of the choice of the 
element in $\hat \P$, we conclude that $\hat \P$ is a principal bundle  over $\hat E$, 
with structure group isomorphic to the group 
of matrices 
$$G = \left\{ \ \left( \matrix 
\smallmatrix 1 & 0  & 0 & 0 \\
B_1  & C_1 & 0 & 0  \\
\bar B_1 & 0 & \bar C_1 & 0\\
A_1 & \frac{i}{2} \bar B_1 C_1  &  - \frac{i}{2} B_1 \bar C_1 & 1  \endsmallmatrix & 0 \\
0 & 
\smallmatrix 1 & 0  & 0 & 0 \\
 B_2 & C_2 & 0  & 0 \\
\bar B_2  & 0 & \bar C_2 & 0\\
A_2 & \frac{i}{2} \bar B_2  C_2 &  - \frac{i}{2} B_2 \bar C_2 & 1  \endsmallmatrix \endmatrix \right)\ 
\smallmatrix B_i, C_i\in \C,\  A_i \in \R,\  |C_i| = 1\endsmallmatrix \ ,\ \right\}
\tag 5.16$$
Note that $G$ is isomorphic to a subgroup of the isotropy $\HQ$, where $\Q$ is the hyperbolic 
quadric. \par
\bigskip
\subsubhead 5.4 Conclusion
\endsubsubhead
\medskip
The following Theorems are the main result of this section.\par
\medskip 
\proclaim{Theorem 5.2} Let $M$ be  elliptic or hyperbolic, 
 $G$  the group described in (5.11) and (5.16), respectively, and 
$\hat G$ the  
structure group of $\hat \pi: \hat E \to M$. Then the group $\hat G
\ltimes G$ has a natural right action  on $\hat \P$, determined by the
right action of $\hat G$ given 
in Proposition 5.1 and the right action of $G$ described in \S 5.2 and \S 5.3, and
$\pi: \hat \P \to M$ is a principal bundle over $M$ with 
structure group $\hat G \ltimes G$.\par
Moreover, $\hat G \ltimes G \simeq \HQ \times \Z_2$ if $M$ is elliptic 
or $\hat G \ltimes G \simeq \HQ \times \Z_2 \times \Z_2$ if $M$ is hyperbolic, where  
$\HQ$ is the stability subgroup of the osculating 
quadric  $\Q$ at the origin. 
\endproclaim
\demo{Proof} The group $\HQ$ is described e. g. in \cite{16}. The existence of the
isomorphisms $\hat G \ltimes G \simeq \HQ \times \Z_2$
or $\hat G \ltimes G \simeq \HQ \times \Z_2 \times \Z_2$ 
follows by  comparison of 
the  groups. Finally, the fact that $\hat G \ltimes G$  acts transitively 
on the fibers of $P$ is a direct consequence of 
the definitions of the  actions of $\hat G$ and $G$.\qed
\enddemo
\bigskip
\proclaim{Theorem 5.3} The Chern-Moser bundle $\P$ of an elliptic or
hyperbolic manifold $M$ is a union of connected components of $\hat \P$. 
In particular, it is a principal bundle over $M$ with structure 
group isomorphic to $\HQ$.
\endproclaim
\demo{Proof} For any point $\theta\in \hat E$ and any $u = (e_1, \dots, e_8)\in \hat \P|_\theta$, 
consider the subspace $\tilde \D_{\theta,u} = \ker(e^7) \cap \ker(e^8)$. The
projection
$\hat \pi_*|_{T_\theta\hat E}: T_\theta\hat E \to T_{\hat \pi(\theta)} M$ induces a linear isomorphism between $\tilde \D_{\theta,u}$ and 
$\D_{\hat \pi(\theta)}$. We may consider the induced complex structure $\tilde J_{\theta,u} : 
\tilde \D_{\theta, u} \to \tilde \D_{\theta, u}$ defined by 
$$\tilde J_{\theta,u} = (\hat \pi_*|_{\tilde \D_{\theta,u}})^{-1} \circ J \circ 
\hat \pi_*|_{\tilde \D_{\theta,u}}\ .$$
So, a frame $u = (e_1, \dots, e_8)\in \hat \P$ belongs to $\P$ 
if and only if 
$$e^4(\tilde J_{\pi_o(u),u}(e_3)) = +1\ ,\qquad e^6(\tilde J_{\pi_o(u),u}(e_5)) = +1\ .$$
Since the functions on the left hand sides are continuous and take values $\pm 1$, 
they are  constant on any connected component of 
$\hat \P$ and this immediately implies the first claim. \par
The second claim follows from the fact that, for any 
$x\in M$, the intersection  $\hat \P_x \cap \P$
coincides with the orbit of $\HQ \subset \hat G \ltimes G$. \qed
\enddemo
\bigskip
\bigskip
\subhead 6. A canonical Cartan connection 
on the Chern-Moser bundle
\endsubhead
\bigskip
\definition{Definition 6.1} Let $P$ be an $H$-principal bundle over 
$M$ and assume that there exists a Lie algebra $\g$, which properly contains $\h = Lie(H)$
and a representation $\operatorname{Ad}: H \to \operatorname{Aut}(\g)$ which extends 
the adjoint representation of $\h$ on $\g$.\par
A {\it Cartan connection on $P$ with model $(\g,\operatorname{Ad}(H))$\/}
is a $\g$-valued 1-form $\psi: TP \to \g$ such that (here, for any $h\in H$, $R_h$ is the right action of $h$ on 
$P$): 
\roster
\item"i)" for any $u\in P$, $\psi_u: T_uP \to \g$ is a linear isomorphism; 
\item"ii)" $\psi(A^*) = A$ for any $A \in \goth h = Lie(H)$ (here $A^*$
is the fundamental vector field  associated with $A$, i.e. 
the vector field on $P$, whose flow is equal to the 1-parameter family of diffeomorphisms
$R_{\exp(tA)}$)
\item"iii)" $\psi$ is $H$-invariant, i.e. for any $h\in H$
$$R_h^*\psi = \operatorname{Ad}_{h^{-1}} \circ \psi\ .\tag6.1$$ 
\endroster
\enddefinition
As before, for a given elliptic or hyperbolic manifold $M$, 
we will always denote  by $\GQ$ and $\HQ$ 
the  group of automorphisms of the osculating quadric $\Q$  and the stability 
subgroup at the origin, respectively.   
Finally, we set $\gQ = Lie(\GQ)$ and $\hQ = Lie(\HQ)$. \par 
\medskip
By standard facts on elliptic or hyperbolic quadrics (see e.g. \cite{16}), 
it is known that the Lie algebra $\gQ = Lie(\GQ)$ is a semisimple  Lie algebra 
isomorphic to $\goth{su}_{2,1} \oplus \goth{su}_{2,1}$ (hyperbolic case) or 
to $\goth{sl}_3(\C)$ (elliptic case). It also admits 
a graded decomposition  of the form
$$\gQ = \gQ^{-2} + \gQ^{-1} + \gQ^0 + \gQ^1 + \gQ^2\ ,$$
with $[\gQ^i, \gQ^j] \subset \gQ^{i+j}$ such that 
 $\hQ = \gQ^0 + \gQ^1 + \gQ^2$ and so that the following property holds: 
if for any $X\in \gQ$ we denote by 
$\hat X$ the corresponding infinitesimal transformation on 
$\Q \subset \GQ/\HQ$, then the map 
$$\imath: \gQ \to T_0 \Q\ ,\qquad X \overset{\imath}\to\mapsto \hat X_{e\HQ}\tag 6.2$$
induces  an isomorphism between
$\gQ^{-1}$ and the holomorphic tangent space $\D_0$ and 
an isomorphism between $\gQ^{-2} + \gQ^{-1}$ and $T_0 \Q$.\par 
Our aim is to construct explicitly  a  Cartan connection on $\P$ with model 
$(\gQ$, $\operatorname{Ad}(\HQ))$,
which is invariant under any CR transformation of $M$. Any  
connection which satisfies such   property of invariance will be called  {\it canonical\/}. \par
\medskip
From now on, we will denote by $(\epsilon_i,  V_{A,k})$ the special 
basis for $\gQ$ listed 
in Appendix. Notice that such a basis is so that: 
\roster
\item"i)" $(\epsilon_1, \epsilon_2)$ is a basis for $\gQ^{-2}$, 
$(\epsilon_3, \dots, \epsilon_6)$ is a basis for the subspace 
$\gQ^{-1} \subset \gQ$  and 
$$J \imath(\epsilon_3) = \imath(\epsilon_4)\ ,\quad 
J \imath(\epsilon_5) = \imath(\epsilon_6)\ $$
(here $\imath$ is the map defined in (6.2) and $J$ is the complex structure of $T_0 \Q$); 
\item"ii)"  for any $k = 0, 1, 2$, the elements  $(V_{A,k})$ form a basis for the subspace 
$\gQ^{k} \subset \gQ$; 
\item"iii)"  the  Lie brackets between the elements $\epsilon_i$ are as follows:
\endroster 
\medskip
\moveright 0.5 cm
\vbox{\offinterlineskip
\halign {\strut\vrule\hfil\ $#$\ \hfil
 &\vrule\ \vrule\hfil\ $#$\ \hfil
\vrule\cr
\noalign{\hrule}
\text{Elliptic}
& \overset{\phantom{AAA}}\to{
\underset{\phantom{AAA}}\to{
\matrix [\epsilon_3, \epsilon_5] = - \epsilon_1\ ,\quad [\epsilon_4, \epsilon_6] = - \epsilon_1\ ,\quad
[\epsilon_3, \epsilon_6] = - \epsilon_2\ ,\quad [\epsilon_4, \epsilon_5] = \epsilon_2\\
\phantom{A}\\
[\epsilon_3, \epsilon_4]Ê= [\epsilon_5, \epsilon_6] = 0
\endmatrix
}}
\cr \noalign{\hrule}
\text{Hyperbolic}
&
\overset{\phantom{AAA}}\to{
\underset{\phantom{AAA}}\to{
\matrix
[\epsilon_3, \epsilon_4] = - \epsilon_1\ ,\quad [\epsilon_5,\epsilon_6] = - \epsilon_2\ ,\\
\phantom{AA}\\
[\epsilon_3, \epsilon_5] = [\epsilon_3, \epsilon_6] = [\epsilon_4, \epsilon_5] =
[\epsilon_4, \epsilon_6] = 0
\endmatrix }}
\cr \noalign{\hrule}
}}
\smallskip
\centerline{\bf Table 5}
\medskip
A given $\gQ$-valued 1-form $\psi$ on $\P$ can be always written  as
$$\psi = \sum_i \epsilon_i \psi^{\epsilon_i} + \sum V_{A,k} \psi^{V_{A,k}}$$
where $\psi^{\epsilon_i}$ and $\psi^{V_{A,k}}$ denote some suitable $\R$-valued 1-forms on $\P$.\par
Observe that a  $\gQ$-valued 1-form $\psi$ satisfies (i) of Definition 6.1
if and only if for any $u\in \P$ 
the 1-forms $\psi^{\epsilon_i}|_u$ and 
$\psi^{V_{A,k}}|_u$ are
a basis for $T^*_u\P$. If this occurs, we have a natural injective linear homomorphism
 between $\gQ$ and the vectors fields on $\P$, namely the correspondence between 
any element $X\in \gQ$ and the unique vector field $\hat X$ such that $\psi_u(\hat X) = X$ at 
any $u\in \P$. Such a vector field $\hat X$ will be called 
{\it fundamental vector fields associated with $X \in \gQ$ by means of the $\gQ$-valued 
1-form $\psi$\/}.\par
\medskip
Then next Lemma gives a characterization of the Cartan connections amongst the $\gQ$-valued 
1-forms which satisfy Definition 6.1 (i). \par 
\medskip
\proclaim{Lemma 6.2} Let $(\psi^{\epsilon_i}, \psi^{V_{A,k}})$ be a set of 
$\R$-valued 1-forms on $\P$, which are linearly independent at all points of $\P$, and let
$\psi = \sum_i \epsilon_i \psi^{\epsilon_i} + \sum V_{A,k} \psi^{V_{A,k}}$. For any $X\in \gQ$, 
denote also by 
$\hat X$ the associated fundamental vector field, by means of $\psi$. Then $\psi$
is a Cartan connection modelled on $(\gQ, \operatorname{Ad}(\HQ))$
if and only if 
\roster
\item the vector fields $\hat V_{A,k}$ coincide with the 
fundamental vector fields $V^*_{A,k}$ associated with the 
elements $V_{A, k} \in \hQ \subset \gQ$, by means of the right action of $\HQ$ on $\P$; 
\item for any vector field $\hat \epsilon_i$
and  any element $V_{B,k} \in \hQ$ 
$$[V_{B,k}, \epsilon_i] = - d\psi(V^*_{B,k}, \hat \epsilon_i) = \psi([V^*_{B,k}, \hat \epsilon_i])\ .\tag 6.3$$
\item $\psi$ is invariant under the element $\gsym \in \HQ \subset \hat G \ltimes G$ defined as follows: if 
$M$ is elliptic, then $\gsym = g_1 \cdot g_2$, where $g_1 = \left(\smallmatrix 1 & 0\\ 0 & -1 
\endsmallmatrix\right) \in \hat G$ and $g_2 \in G$ is the element 
which exchanges $\bar E^1$ with $E^2$; if $M$ is hyperbolic, then
$\gsym = \left(\smallmatrix 0 & 1\\ 1 & 0
\endsmallmatrix\right) \in \hat G$.
\endroster
\endproclaim
\demo{Proof} By construction, $\psi$ satisfies i) of Definition 6.1. 
Moreover, if (1) holds,  then $\psi$ satisfies ii) of Definition 6.1 and,
 for any $V_{A, k} , V_{B,\ell} \in \hQ$ we have that
$$[V_{B,\ell}, \psi(\hat V_{A, k})]  + (\Cal L_{V^*_{B,\ell}} \psi)(\hat V_{A, k}) = 
[V_{B,\ell}, V_{A,k}] - \psi([V^*_{B,\ell}, V^*_{A,k}]) = $$
$$ = [V_{B,\ell}, V_{A,k}] - \psi([V_{B,\ell}, V_{A,k}]^*) = 0\ .
\tag 6.4$$
It follows that, if $\psi$ satisfies also (6.3), 
$$[V_{B,k}, \psi(\hat\epsilon_i)] + (\Cal L_{V^*_{B,\ell}} \psi)(\hat \epsilon) 
= [V_{B,k}, \epsilon_i]  +  d\psi(V^*_{B,k}, \hat \epsilon_i) = 0\ .$$
So,  for any element $\hat X \in 
\{\hat \epsilon_i, \hat V_{A,k}\}$, we have that 
$[V_{B,k}, \psi(\hat X)]  + (\Cal L_{V^*_{B,k}} \psi)(\hat X)  = 0$.
Since the vectors $\{\hat \epsilon_i, \hat V_{A,k}\}$ span $T_u\P$ at any point and, in the elliptic 
and hyperbolic case, the connected component of the identity $\HQ^o$
coincides with $\exp(\hQ)$, we conclude that  $\psi$ is 
$\HQ^o$-invariant. Finally, since $\HQ = \HQ^o \times \{e, \gsym\}$, if also (3) is satisfied, 
then $\psi$ is invariant under the entire group $\HQ$ and hence it is a Cartan connection.\par
 The necessity of conditions (1) - (3)  follows from the definitions. 
\qed
\enddemo
\medskip
\remark{Remark 6.3} 
Note that condition (6.3) can be written also as
$$[V^*_{B,k}, \hat \epsilon_i] = \widehat{[V_{B,k}, \epsilon_i]}\ .\tag 6.3'$$
Observe also that, by the proof of Proposition 5.1, if $M$ is hyperbolic, the element 
$\gsym$ maps any adapted coframe $(e^1, \dots, e^8)$
into the adapted coframe 
$(e'{}^1, \dots, e'{}^8)$ defined by 
$$e'{}^1 = R_{\gsym}^*(e^2)\ ,\ e'{}^2 = R_{\gsym}^*(e^1)\ ,$$
$$e'{}^3 = R_{\gsym}^*(e^5)\ ,\ e'{}^4 = R_{\gsym}^*(e^6)
\ ,\ e'{}^5 = R_{\gsym}^*(e^3)\ ,\ e'{}^6 = R_{\gsym}^*(e^4)\ .$$
Similarly, by the same Proposition 5.1, if $M$ is elliptic, 
the element 
$\gsym$ maps any adapted coframe 
$(e^1, \dots, e^8)$
into the adapted coframe 
$(e'{}^1, \dots, e'{}^8)$ defined by 
$$E'{}^0  = R_{\gsym}^*(\bar E^0) \ ,\ E'{}^1  = R_{\gsym}^*(E^2) \ ,
\ E'{}^2 = R_{\gsym}^*(E^1)\ ,$$
where, as before,  $E^0 = e^1 + i e^2$, $E^1 = e^3 + i e^4$ and $E^2 = e^5 + i e^6$. 
\endremark
\medskip
The next lemma shows that,  locally, a Cartan connection always exists
and that some additional useful properties can be always assumed.\par
\medskip
\proclaim{Lemma 6.4} Let  $\Cal U \subset M$ be an open subset 
such that $\pi^{-1}(\Cal U) \subset \P$ is trivializable.  Then there exists a Cartan connection $\psi$ on 
$\pi^{-1}(\Cal U)$
modelled on $(\gQ, \operatorname{Ad}(\HQ))$, such that 
$$\psi^{\epsilon_i} = \omega^i\ ,\tag 6.5$$
where the $\omega^i$'s are the components of the 
tautological 1-form $\omega$ of $\P$.
\endproclaim
\demo{Proof} It is known (see e.g. Ch. 5 in \cite{17}) that any Cartan connection on 
$\pi^{-1}(\Cal U)$ can be constructed as follows. Let $\mu: \pi^{-1}(\Cal U) \to \Cal U \times \HQ$
be a trivializing map and let $\delta_x$  be a  family of 
linear maps 
$\delta_x : T_x \Cal U \to \gQ$,
depending smoothly on the points $x\in \Cal U$ and 
such that the compositions  $p \circ \delta_x$ with the projection $p: \gQ \to \gQ^{-2} + \gQ^{-1}$ are 
linear isomorphisms. If we denote by $\pi_2: \Cal U \times 
\HQ \to \HQ$ the natural projection onto the 
second factor and we denote by
$\omega_{\HQ}$  the Maurer-Cartan form of $\HQ$, then the 1-form
$$\psi_{\mu^{-1}(x,h)} = \mu^*\left(
\operatorname{Ad}_{h^{-1}} \circ \delta_x + \pi_2^*\omega_{\HQ} \right)\tag 6.6$$
is a Cartan connection. Moreover,  any Cartan connection that is
modelled on $(\gQ, \operatorname{Ad}(\HQ))$
on   $\pi^{-1}(\Cal U)$ is of the above form. \par
Consider now the family of linear maps 
$$\delta_x : T_x \Cal U \to \gQ^{-2} + \gQ^{-1}\  ,\qquad \delta_x(X) = \sum \epsilon_i \omega^i_{\mu^{-1}(x,e)}
(\mu^{-1}_*(X))$$
and the corresponding Cartan connection $\psi$ defined by (6.6). By construction,
(6.5) holds  at all 
points of the form $\mu^{-1}(x,e)$. By the $\HQ$-invariance of 
$\psi$, the properties of the vectors $\epsilon_i$ and the transformation rules of the tautological 
1-form 
$\omega$, it follows that the identity (6.5)  is  satisfied at any  
point of $\pi^{-1}(\Cal U)$. 
\qed 
\enddemo
\bigskip
A (local) Cartan connection which satisfies (6.5) will be called {\it good\/}.\par
Consider a (local) good Cartan connection $\psi$
and let $\hat \epsilon_i$ and 
$\hat V_{A, k}$ be the associated fundamental vector fields. 
It is not hard to check  that, if $\psi'$ is a new 
$\gQ$-valued 1-form, which   satisfies (6.5) and  conditions (i) and (ii) of Definition 
6.1, then there  exist some smooth $\R$-valued functions $S_i^{A,k}$ such that 
the fundamental vector fields $\hat \epsilon'_i$, 
$\hat V'_{A,k}$, determined by $\psi'$,  and the components $\psi'{}^{\epsilon_i}$, 
$\psi'{}^{V_{A,k}}$ of $\psi'$ are as follows:
$$\hat \epsilon'_i = \hat \epsilon_i + \sum_{A,k} S_i^{A,k} \hat V_{A,k}\ ,\qquad
\hat V'_{A, k}= \hat V_{A,k} = V^*_{A,k}\ ,\tag 6.7$$
$$\psi'{}^{\epsilon_i} = \omega^i\ ,\qquad 
\psi'{}^{V_{A,k}} = \psi^{V_{A,k}} - \sum_i S_i^{A,k}\omega^i\ .\tag 6.8$$
On the other hand, by Lemma 6.2, this new $\gQ$-valued 1-form  $\psi'$ is a 
Cartan connection if and only if it is invariant 
under the element $\gsym \in \HQ$ and equation
(6.3) is satisfied. This last condition means  that 
for any element $V_{B,\ell}$ and any $\epsilon_i$ the following has to be satisfied:  
$$[V_{B,\ell}, \epsilon_i] = 
\sum_j \epsilon_j\omega^j([V^*_{B,\ell}, \hat \epsilon'_i]) + \sum_{C, m} V_{C,m} \psi'{}^{V_{C,m}}(
[V^*_{B,\ell}, \hat \epsilon'_i]) = $$
$$ = \sum_j \epsilon_j\omega^j([V^*_{B,\ell}, \hat \epsilon_i]) + 
\sum_{j, A,k} \epsilon_j S_i^{A,k} \omega^j([V^*_{B,\ell}, V^*_{A,k}]) +
\sum_{j, A,k} \epsilon_j V^*_{B,\ell}\left(S_i^{A,k}\right) \omega^j( V^*_{A,k}) + $$
$$ + \sum_{C, m} V_{C,m} \left(\psi^{V_{C,m}} - \sum_j S^{C, m}_j \omega^j\right)(
[V^*_{B,\ell}, \hat \epsilon_i + \sum_{A,k} S^{A,k}_i V_{A,k}^*]) = $$
$$ =\sum_j \epsilon_j\omega^j([V^*_{B,\ell}, \hat \epsilon_i]) + \sum_{C, m} V_{C,m} \psi^{V_{C,m}}
([V^*_{B,\ell}, \hat \epsilon_i]) + $$
$$ + \sum_{A,k, C, m} V_{C,m} 
V^*_{B,\ell}(S^{A,k}_i) \psi^{V_{C,m}}(V^*_{A,k}) + $$
$$ + \sum_{A,k,C, m} S^{A,k}_i V_{C,m} \psi^{V_{C,m}}([V^*_{B,\ell},V_{A,k}^*])
- \sum_{j, C, m} V_{C,m} S^{C, m}_j \omega^j([V^*_{B,\ell}, \hat \epsilon_i])\ .\tag 6.9$$
Formula 
(6.9) simplifies considerably  if we 
recall that  $\psi$ is a Cartan connection (and hence satisfies (6.3)). In this way, we obtain that
$\psi'$ is a Cartan connection if and only if it is $\gsym$-invariant and 
the functions $S^{A,k}_i$ satisfy, for any $V^*_{B,\ell}$
$$V^*_{B,\ell}(S^{A,k}_i) = \sum_j  S^{A, k}_j \omega^j([V^*_{B,\ell}, \hat \epsilon_i]) - 
\sum_{C,m} S^{C,m}_i \psi^{V_{A,k}}([V^*_{B,\ell},  V_{C,m}^*]) 
\ .\tag 6.10$$
In the next two subsection, we will show that, for any  good Cartan connection $\psi$,
defined on a trivializable subset   $\pi^{-1}(\Cal U) \subset \P$,
there exists a  {\it unique\/} choice for a 
1-form $\psi'$ presented as in (6.7) and (6.8)  such that: 
\roster
\item"a)" certain conditions 
on the values $d\omega^i(\hat \epsilon'_j, \hat \epsilon'_k)$
and $d\psi'{}^{V_{A,k}}(\hat \epsilon'_j, \hat \epsilon'_k)$ are satisfied;
\item"b)" $\psi'$ is $\gsym$-invariant 
and satisfies (6.10) (hence it is a good Cartan connection). 
\endroster
The existence and uniqueness of such   modification $\psi'$ for any 
given  good Cartan connection implies that it has to coincide on the 
overlaps of two trivializable sets  $\pi^{-1}(\Cal U)$, $\pi^{-1}(\Cal U')$, 
even if 
we started from   two distinct good Cartan connections $\psi_1$ and $\psi_2$. 
For this reason all modifications $\psi'$  can be patched together to defined a unique 
global good Cartan connection $\psi_{CM}$ on $\P$, which is necessarily canonical. \par
\medskip
Finally, let us explain how it is possible 
to check if the canonical connection we obtain is the same (or, more 
precisely, equivalent) to the connection $\omega_M$ of Theorem 1.1.
This can be done 
by means of  Theorem 2.7 in 
\cite{20}, where the following necessary and sufficient conditions 
for the existence of an  isomorphism between $\psi_{CM}$ and $\omega_M$ are given.
Let us denote by 
$(\varepsilon_A)$, $1 \leq A \leq 6$, a basis for the subspace 
$\m \= \gQ^{-2} + \gQ^{-1} \subset \gQ$ and let $(\varepsilon^A)$ a corresponding basis 
for $\gQ^1 + \gQ^2 \subset \gQ$, which is dual w.r.t. the Cartan-Killing form $\Cal B$ of $\gQ$, i.e. 
such that 
$\Cal B(\varepsilon^A, \varepsilon_B) = \delta^A_B$. Then,  in Theorem 2.7 in \cite{20}
it is given a necessary and sufficient condition,  which may be rephrased for 
 good Cartan connections saying that
 the following  equations are satisfied for any $e_B$ (to obtain the following  
expression from the original statement in \cite{20}, we used the fact that
for any 
$X\in \gQ$ and any $\hat e^A$, 
$[\hat e_A, \psi_{CM}(\hat X)] = - {\Cal L}_{\hat e^A}\psi_{CM}(\hat X)$  - see proof of Lemma 6.2):
$$\sum_A\hat\varepsilon^A\left(d\psi_{CM}(\hat \varepsilon_A, \hat \varepsilon_B) \right) -
\sum_A \frac{1}{2} d\psi_{CM}(\widehat{[\varepsilon^A, \varepsilon_B]_\m}, \hat \varepsilon_A ) = 
\sum_A \frac{1}{2}[[\varepsilon^A, \varepsilon_B]_\m, \varepsilon_A] \ ,\tag 6.11$$
where $[\varepsilon^A, \varepsilon_B]_\m$ denotes the natural projection of $[\varepsilon^A, \varepsilon_B]$
into the subspace $\m = \gQ^{-2} + \gQ^{-1}$.\par
If one consider a basis $(\varepsilon_A)$, $1 \leq A \leq 6$, given by the first 
six elements of 
the  special basis for $\gQ$ described in the Appendix, it is possible to realize that 
the corresponding  dual basis $(\varepsilon^A)$ is given (up to  factors) by 
the last six elements of the same basis for $\gQ$. Then, using (6.12) and (6.13) below, 
it is possible to write down explicitly  all components of the $\gQ$-valued 1-form on the left hand 
side of (6.11) and determine a set of
 conditions which is equivalent to (6.11). We will see 
that, among them, there are some of the 
which are not satisfied by a generic connection
defined by the conditions mentioned in (a). \par
\medskip
Let us now proceed with the 
construction of the canonical Cartan connection $\psi_{CM}$, 
following the steps (a) and (b) described above. \par
Before going into the details of such a construction, we need
the following technical fact, whose proof is just an application of definitions 
and of (6.3). In the following, most of the 
functions  $S^{A,k}_i$ will coincide with  functions  which are 
linear combinations of the functions $d\omega^k(\hat \epsilon_i, 
\hat \epsilon_j)$ and $d\psi^{V_{C,0}}(\hat \epsilon_i, \hat \epsilon_j)$. In order to check if 
(6.10) holds we need to evaluate the directional derivatives of such functions 
and they are given by the following expressions (we are assuming that $\psi$ is 
a Cartan connection): 
$$V_{B,m}^*\left(d\omega^k(\hat\epsilon_i, \hat \epsilon_j)\right) = 
d^2 \omega^k(V_{B,m}^*, \hat\epsilon_i, \hat \epsilon_j) - 
\hat \epsilon_i\left(d\omega^k(\hat\epsilon_j, V_{B,m}^*)\right)  - 
\hat \epsilon_j\left(d\omega^k(V_{B,m}^*, \hat\epsilon_i)\right) + $$
$$ + 
d\omega^k([V_{B,m}^*, \hat \epsilon_i], \hat \epsilon_j) + 
d\omega^k([\hat \epsilon_i, \hat \epsilon_j],V_{B,m}^*) + 
d\omega^k([\hat \epsilon_j, V_{B,m}^* ], \hat \epsilon_i) \overset{(6.3)}\to= $$
$$ = d\omega^k([V_{B,m}^*, \hat \epsilon_i], \hat \epsilon_j) + 
d\omega^k([\hat \epsilon_i, \hat \epsilon_j],V_{B,m}^*) + 
d\omega^k([\hat \epsilon_j, V_{B,m}^* ], \hat \epsilon_i) \overset{(6.3)}\to=$$
$$ = d\omega^k(\widehat{[V_{B,m}, \epsilon_i]}, \hat \epsilon_j) - 
d\omega^k(\widehat{[V_{B,m},  \epsilon_j ]}, \hat \epsilon_i) + 
\sum_\ell\omega^\ell([\hat \epsilon_i, \hat \epsilon_j])
d\omega^k(\hat \epsilon_\ell, V_{B,m}^*) + $$
$$ + 
\sum_{C,\ell} \psi^{V_{C,\ell}}([\hat \epsilon_i, \hat \epsilon_j]) d\omega^k(V_{C,\ell}^*, V_{B,m}^*) = $$
$$ = d\omega^k(\widehat{[V_{B,m}, \epsilon_i]}, \hat \epsilon_j) + 
d\omega^k(\hat \epsilon_i, \widehat{[V_{B,m},  \epsilon_j ]}) - \sum_{\ell} d\omega^\ell(\hat \epsilon_i, \hat \epsilon_j)
\omega^k(\widehat{[V_{B,m}, \epsilon_\ell]})\ ,\tag 6.12$$
$$V_{B,m}^*\left(d\psi^{V_{C,\ell}}(\hat\epsilon_i, \hat \epsilon_j)\right) = $$
$$= d\psi^{V_{C,\ell}}([V_{B,m}^*, \hat \epsilon_i], \hat \epsilon_j) + 
d\psi^{V_{C,\ell}}([\hat \epsilon_i, \hat \epsilon_j],V_{B,m}^*) + 
d\psi^{V_{C,\ell}}([\hat \epsilon_j, V_{B,m}^* ], \hat \epsilon_i)  = $$
$$= d\psi^{V_{C,\ell}}(\widehat{[V_{B,m}, \epsilon_i]}, \hat \epsilon_j) + 
d\psi^{V_{C,\ell}}(\hat \epsilon_i, \widehat{[V_{B,m},  \epsilon_j ]}) - $$
$$ - 
\sum_{n} d\omega^n(\hat \epsilon_i, \hat \epsilon_j)
\psi^{V_{C,\ell}}(\widehat{[V_{B,m}, \epsilon_n]}) - 
\sum_{A,k} d\psi^{V_{A,k}}(\hat \epsilon_i, \hat \epsilon_j)
\psi^{V_{C,\ell}}(\widehat{[V_{B,m}, V_{A,k}]})\ .
\tag 6.13$$
\bigskip
\subsubhead 6.1 Construction of a canonical Cartan connection on a hyperbolic manifold
\endsubsubhead
\medskip
Assume that $M$ is hyperbolic. In this case $\gQ \simeq \su_{2,1} \oplus \su_{2,1}$
and  the special basis   $(\epsilon_i, V_{A,k})$  is given in Appendix.
Notice also that,  using the same arguments which brought to  (4.5),
any good Cartan connection $\psi$ is so that
$$d\omega^1 =  \omega^3 \wedge \omega^4 + S \omega^5 \wedge \omega^2 +  T \omega^6 \wedge \omega^2 
+ 2\psi^{V_{2,1,0}}\wedge\omega^1\ \  
$$
$$+\ \text{linear combinations of }\ \ \left\{\ \omega^1\wedge \omega^2\ ,
\ \omega^1 \wedge\omega^j\ ,\ 3\leq j\leq 6\ \right\}\tag 6.14$$
$$d\omega^2 =  \omega^5 \wedge \omega^6 + S' \omega^3 \wedge \omega^1 +  T' \omega^4 \wedge \omega^1 + 
2\psi^{V_{2,2,0}}\wedge\omega^2\ \  
$$
$$+\  \text{linear combinations of }\ \ \left\{
\ \omega^1\wedge \omega^2\ ,\ \omega^2 \wedge\omega^j\ ,\  3\leq j\leq 6\ \right\}\tag 6.15$$
for some suitable functions $S,S', T, T'$. Moreover, from the vanishing of 
$$d^2\omega^1({\hat \epsilon}_{4},{\hat \epsilon}_{5},{\hat \epsilon}_{6}) =
d^2\omega^1(\hat \epsilon_3, \hat \epsilon_5, \hat \epsilon_6) = 
d^2\omega^2({\hat \epsilon}_{3},{\hat \epsilon}_{4},{\hat \epsilon}_{6}) =  
d^2\omega^2({\hat \epsilon}_{3},{\hat \epsilon}_{4},{\hat \epsilon}_{5}) = 0$$ 
and using (6.3), 
one can check that for any good Cartan connection $\psi$ the following identities hold:
$$d\omega^3(\hat \epsilon_5, \hat \epsilon_6) = d\omega^4(\hat \epsilon_5, \hat \epsilon_6)  = 
d\omega^5(\hat \epsilon_3, \hat \epsilon_4) = d\omega^6(\hat \epsilon_3, \hat \epsilon_4) = 0\ .\tag 6.16$$
\par
Now we will proceed with the construction of a canonical Cartan connection, which is based 
on a sequence of technical lemmata. 
\medskip
\proclaim{Lemma 6.5} On any trivializable open subset $\pi^{-1}(\Cal U) \subset \P$, 
there exists a good Cartan connections $\psi$, which 
satisfies the following conditions  for any  $a = 1,2$, $j = 1,3,4$ and $k = 2,5,6$:
$$d\omega^1(\hat \epsilon_1, \hat \epsilon_3)   = 
d\omega^1(\hat \epsilon_1, \hat \epsilon_4) = 
d\omega^2(\hat \epsilon_2, \hat \epsilon_5) =
d\omega^2(\hat \epsilon_2, \hat \epsilon_6) = 
0\ ,\tag 6.17$$
$$d\omega^3(\hat \epsilon_3, \hat \epsilon_4) =  d\omega^4(\hat \epsilon_3, \hat \epsilon_4) = 
d\omega^5(\hat \epsilon_5, \hat \epsilon_6)  = 
d\omega^6(\hat \epsilon_5, \hat \epsilon_6) =  
0\ ,\tag 6.18$$
$$d\omega^6(\hat \epsilon_j, \hat \epsilon_5) - d\omega^5(\hat \epsilon_j, \hat \epsilon_6)
= d\omega^5(\hat \epsilon_j, \hat \epsilon_5) + d\omega^6(\hat \epsilon_j, \hat \epsilon_6) + 
d\omega^2(\hat \epsilon_j, \hat \epsilon_2)= 0\ ,\tag 6.19$$
$$d\omega^4(\hat \epsilon_k, \hat \epsilon_3) - d\omega^3(\hat \epsilon_k, \hat \epsilon_4)
= d\omega^3(\hat \epsilon_k, \hat \epsilon_3) + d\omega^4(\hat \epsilon_k, \hat \epsilon_4) + 
d\omega^1(\hat \epsilon_k, \hat \epsilon_1)= 0\ ,\tag 6.20$$
$$d\omega^3(\hat \epsilon_3, \hat \epsilon_1) = d\omega^4(\hat \epsilon_4, \hat \epsilon_1) =
d\omega^3(\hat \epsilon_4, \hat \epsilon_1) = d\omega^4(\hat \epsilon_3, \hat \epsilon_1) = 0\ ,\tag 6.21$$
$$d\omega^5(\hat \epsilon_5, \hat \epsilon_2) = d\omega^6(\hat \epsilon_6, \hat \epsilon_2) = 
d\omega^5(\hat \epsilon_6, \hat \epsilon_2) = d\omega^6(\hat \epsilon_5, \hat \epsilon_2) = 0\ .
\tag  6.22$$
Moreover, if $\psi'$ is another good Cartan connection with the same properties and expressed in terms of 
$\psi$ as in (6.7), (6.8),  
then the  
functions
$S^{a,b,0}_J$, $S^{2,2,0}_1$, 
$S^{2,1,0}_2$,   with $a, b, = 1,2$ and $3\leq J\leq 6$, are vanishing, while the others 
satisfy the following linear relations: 
$$S^{1,1,0}_1 =  S^{2,1,1}_4 = - S^{1,1,1}_3\ ,\quad S^{2,1,0}_1 = S^{2,1,1}_3 = S^{1,1,1}_4\ ,\tag 6.23$$
$$S^{1,2,0}_2 = S^{2,2,1}_6 = - S^{1,2,1}_5\ ,\quad S^{2,2,0}_2 = S^{2,2,1}_5 = S^{1,2,1}_6\ .\tag 6.24$$
Finally, for any good Cartan connection $\psi$, which satisfies (6.17) - (6.21), the following identities hold
for any $ a = 1,2$, $j = 1,3,4$ and $k = 2,5,6$:
$$ 
d\omega^2(\hat \epsilon_2,\hat \epsilon_j) = d\omega^1(\hat \epsilon_1,\hat \epsilon_k) = 
d\psi^{V_{2,a,0}}(\hat \epsilon_3, \hat \epsilon_4) = d\psi^{V_{2,a,0}}(\hat \epsilon_5, \hat \epsilon_6) =
0\ .
\tag 6.25$$
\endproclaim
\demo{Proof} Consider a good Cartan connection $\psi$ on $\pi^{-1}(\Cal U) \subset \P$
and let $\psi'$ be 
another $\gQ$-valued 1-form defined  as in (6.7) and (6.8) by means of some functions $S^{a,b,k}_i$. 
The 1-form  $\psi'$ satisfies 
$(6.17)_1$ if and only if
$$0 = d\omega^1(\hat \epsilon'_1, \hat \epsilon'_3) = 
- \omega^1([\hat \epsilon_1 + \sum_{b,c,k}S^{b,c,k}_1 V_{b,c,k}^*, \hat \epsilon_3 + 
\sum_{d,e,\ell}S^{d,e,\ell}_3 V_{d,e,\ell}^*]) =$$
$$ = - \omega^1([\hat \epsilon_1 , \hat \epsilon_3])  + S^{2,a,0}_3 \omega^1([V^*_{2,a,0}, \hat \epsilon_1])
= $$
$$ = 
d\omega^1(\hat \epsilon_1, \hat \epsilon_3) - 2 S^{2,1,0}_3\ ,\tag 6.26$$
i.e. 
$S^{2,1,0}_3 = - \frac{1}{2} d\omega^1(\hat\epsilon_3, \hat\epsilon_1)$. 
In a similar way, one can  check that the other conditions in (6.17) and (6.18) are satisfied if and only if 
$$S^{2,1,0}_4 = - \frac{1}{2} d\omega^1(\hat\epsilon_4, \hat\epsilon_1)\ ,\quad
S^{2,2,0}_5 = - \frac{1}{2} d\omega^2(\hat\epsilon_5, \hat\epsilon_2)\ ,
\quad
S^{2,2,0}_6 = - \frac{1}{2} d\omega^2(\hat\epsilon_6, \hat\epsilon_2)\ ,\tag 6.27$$
$$S^{1,1,0}_3 = - d\omega^3(\hat \epsilon_3, \hat \epsilon_4) - \frac{1}{2} d\omega^1(\hat \epsilon_4, \hat \epsilon_1)
\ ,\quad 
S^{1,1,0}_4 =  - d\omega^4(\hat \epsilon_3, \hat \epsilon_4) + \frac{1}{2} d\omega^1(\hat\epsilon_3, \hat\epsilon_1)
\ ,
\tag 6.28$$
$$S^{1,2,0}_j =  
\frac{1}{2}\left(d\omega^6(\hat \epsilon_j, \hat \epsilon_5) - d\omega^5(\hat \epsilon_j, \hat \epsilon_6)  \right)
\ ,\  S^{1,2,0}_j =  
\frac{1}{2}\left(d\omega^4(\hat \epsilon_k, \hat \epsilon_3) - d\omega^3(\hat \epsilon_k, \hat \epsilon_4)  \right)\ ,
\tag 6.29
$$
$$
S^{2,2,0}_j = - 
\frac{1}{4}\left(d\omega^5(\hat \epsilon_j, \hat \epsilon_5) + d\omega^6(\hat \epsilon_j, \hat \epsilon_6) + 
d\omega^2(\hat \epsilon_j, \hat \epsilon_2)  \right)\ ,\tag 6.30$$
$$
S^{2,1,0}_k = - 
\frac{1}{4}\left(d\omega^3(\hat \epsilon_k, \hat \epsilon_3) + d\omega^4(\hat \epsilon_k, \hat \epsilon_4) + 
d\omega^1(\hat \epsilon_k, \hat \epsilon_1)  \right)\tag 6.31$$
for any $j = 1,3,4$ and $ k = 2,5,6$.\par
A  tedious but straightforward check, based on (6.12), shows  that  if we choose the 
 functions $S_J^{A,k}$  so that (6.27)- (6.31) hold and all the remaining
functions so that (6.10) 
is satisfied,  then $\psi'$
satisfies   (6.3) and condition (3) of Lemma 6.2. This means that $\psi'$ is a Cartan connection and
hence that we may assume that $\psi$ satisfies also (6.17) - (6.21) (by replacing $\psi$ with $\psi'$). From 
(6.26) - (6.31), it is clear that any other $\psi'$ with such  properties is defined by functions $S_J^{A,k}$ so that
$S^{a,2,0}_1 = S^{a,1,0}_2 = S^{a,b,0}_J =  0$,  for $a, b = 1,2$ and $3 \leq J \leq  6$.\par
 Let us now look  for  those $\psi'$
which satisfy  also (6.21) and (6.22), we see that this occurs if and only if  the following  holds:
$$S^{2,1,0}_1  - S^{2,1,1}_3 = d\omega^3(\hat \epsilon_3, \hat \epsilon_1) \ ,\
S^{1,1,0}_1  - S^{2,1,1}_4 = d\omega^3(\hat \epsilon_4, \hat \epsilon_1)\ ,\tag 6.32$$
$$S^{1,1,0}_1  + S^{1,1,1}_3 = d\omega^4(\hat \epsilon_1, \hat \epsilon_3) \ ,\
S^{2,1,0}_1  - S^{1,1,1}_4 = d\omega^4(\hat \epsilon_4, \hat \epsilon_1)\ ,\tag 6.33$$
$$S^{2,2,0}_2  - S^{2,2,1}_5 = d\omega^5(\hat \epsilon_5, \hat \epsilon_2) \ ,\
S^{1,2,0}_2  - S^{2,2,1}_6 = d\omega^5(\hat \epsilon_6, \hat \epsilon_2)\ ,\tag 6.34$$
$$S^{1,2,0}_2  + S^{1,2,1}_5 = d\omega^6(\hat \epsilon_2, \hat \epsilon_5) \ ,\
S^{2,2,0}_2  - S^{1,2,1}_6 = d\omega^6(\hat \epsilon_6, \hat \epsilon_2)\ .\tag 6.35$$
It is not hard to realize that the system given by the equations (6.32) -(6.35) has maximal 
rank  and it is therefore solvable at any point. Moreover, with some other straightforward
computations, it can be checked that any set of functions which satisfy (6.32) - (6.35) at the
points of a submanifold  transversal to the fibers  can be smoothly extended 
to functions which satisfy  (6.10) and condition (3) of Lemma 6.2 and which at the 
same time are  solutions of the above given system 
of equations. From this we get the existence of  a Cartan $\psi$ which satisfies also (6.21) and (6.22). 
 The second claim is a direct consequence of (6.26) - (6.35).  The last claim 
can be checked by evaluating $d^2\omega^2(\hat \epsilon_j, \hat \epsilon_5, 
\hat \epsilon_6)$, $d^2\omega^1(\hat \epsilon_k, \hat \epsilon_3, 
\hat \epsilon_4)$, 
$d^2\omega^a(\hat \epsilon_a, \hat \epsilon_3, \hat \epsilon_4)$
and  
$d^2\omega^a(\hat \epsilon_a, \hat \epsilon_5, \hat \epsilon_6)$ with $a=1,2$, $j = 1,3,4$, $k= 2,5,6$,  
and recalling that they  have to vanish because of the well-known triviality of  the operator $d^2$. 
\qed
\enddemo
\medskip
\proclaim{Lemma 6.6} On any trivializable open subset $\pi^{-1}(\Cal U) \subset \P$, 
there exists a  good Cartan connection  $\psi$,
which satisfies (6.17)  - (6.21)  and  the conditions:
$$d\psi^{V_{1,1,0}}(\hat \epsilon_3, \hat \epsilon_4) = 
d\psi^{V_{1,2,0}}(\hat \epsilon_5, \hat \epsilon_6) = 0
\ ,\tag 6.36$$
$$d\psi^{V_{1,1,0}}(\hat \epsilon_3, \hat \epsilon_1) + d\psi^{V_{2,1,0}}(\hat \epsilon_4, \hat \epsilon_1) 
 = d\psi^{V_{2,1,0}}(\hat \epsilon_3, \hat \epsilon_1)  - d\psi^{V_{1,1,0}}(\hat \epsilon_4, \hat \epsilon_1) = 0\ ,
\tag 6.37$$
$$d\psi^{V_{1,2,0}}(\hat \epsilon_5, \hat \epsilon_2) + d\psi^{V_{2,2,0}}(\hat \epsilon_6, \hat \epsilon_2) 
 = d\psi^{V_{2,2,0}}(\hat \epsilon_5, \hat \epsilon_2)  - d\psi^{V_{1,1,0}}(\hat \epsilon_6, \hat \epsilon_2) = 0\ ,
\tag 6.38$$
Moreover, if $\psi'$ is another good Cartan connection with the same properties and expressed in terms of 
$\psi$ as in (6.7), (6.8),  
then all  functions   $S^{A,k}_J$  are vanishing except the functions which appear in $(6.23)_2$, 
$(6.24)_2$ and the functions $S_a^{c,a,1}$, $S_j^{a,2,1}$, $S_k^{a,1,1}$ and
$S_J^{a,b,2}$  with $a = 1,2$, $j= 3,4$, $k = 5,6$ and $1\leq J\leq 6$, among which the 
following relations have to be satisfied:
$$S^{1,1,1}_1 = - \frac{1}{2} S_3^{1,1,2}\ ,\quad S^{2,1,1}_1 = \frac{1}{2} S_4^{1,1,2}\ ,\quad 
S^{1,2,1}_2 = - \frac{1}{2} S_5^{1,2,2}\ ,\quad 
S^{2,2,1}_2 =  \frac{1}{2} S_6^{1,2,2} \ ,\tag 6.39$$
\endproclaim
\demo{Proof} Starting as in the previous lemma, let us
fix a good Cartan connection  $\psi$ on $\pi^{-1}(\Cal U) \subset \P$, which satisfies 
Proposition 6.5,   
and let $\psi'$ be  
an arbitrary $\gQ$-valued 1-form, defined as in (6.8) by means of some functions $S^{a,b,k}_i$, 
which satisfies the same conditions. 
From  the results of the previous lemma on the functions $S^{a,b,k}_i$
and by the same line of arguments used in its proof, 
 we see that $(6.36)_1$  is satisfied if and only if 
$$ 0 = d\psi^{V_{1,1,0}}(\hat \epsilon_3, \hat \epsilon_4) - 
\sum_{i =1}^6 S^{1,1,0}_i d\omega^i(\hat \epsilon_3, \hat \epsilon_4) - 
\frac 3 2 S^{2,1,1}_4 + \frac 3 2 S^{1,1,1}_3 = $$
$$ = d\psi^{V_{1,1,0}}(\hat \epsilon_3, \hat \epsilon_4) -  S^{1,1,0}_1  - 
\frac 3 2 S^{2,1,1}_4 + \frac 3 2 S^{1,1,1}_3 $$
and hence, by (6.23), if and only if 
$$ S^{1,1,0}_1 = S^{2,1,1}_4 = -  S^{1,1,1}_3 = \frac 1 4 d\psi^{V_{1,1,0}}(\hat \epsilon_3, \hat \epsilon_4)\ .\tag 6.40$$
Similarly, the other equality of (6.36) is true if and only if 
$$ S^{1,2,0}_2 = S^{2,2,1}_5 = -  S^{1,1,1}_6 = \frac 1 4 d\psi^{V_{1,2,0}}(\hat \epsilon_5, \hat \epsilon_6)\ .\tag 6.41$$
With similar arguments we get also that (6.37) - (6.42) are satisfied if and only if
$$2 S^{2,1,1}_1 - S_4^{1,1,2} = d\psi^{V_{1,1,0}}(\hat \epsilon_3,\hat \epsilon_1)  + 
d\psi^{V_{2,1,0}}(\hat \epsilon_4,\hat \epsilon_1)\ ,\tag 6.42$$
$$2 S^{1,1,1}_1 + S_3^{1,1,2} = d\psi^{V_{1,1,0}}(\hat \epsilon_4,\hat \epsilon_1) - 
d\psi^{V_{2,1,0}}(\hat \epsilon_3,\hat \epsilon_1)\ ,\tag 6.43$$
$$2 S^{2,2,1}_2 - S_6^{1,2,2} = d\psi^{V_{2,2,0}}(\hat \epsilon_5,\hat \epsilon_2)  + 
d\psi^{V_{1,2,0}}(\hat \epsilon_6,\hat \epsilon_2)\ ,\tag 6.44$$
$$2 S^{1,2,1}_2 + S_5^{1,2,2} = d\psi^{V_{1,2,0}}
(\hat \epsilon_6,\hat \epsilon_2) - d\psi^{V_{2,2,0}}(\hat \epsilon_5,\hat \epsilon_2)\ ,\tag 6.45$$
The proofs of (6.42) - (6.45) are slightly more involved than those in the previous lemma and so, 
for convenience of the
reader, we exhibit the steps which bring to (6.42) (the other conditions are
obtained in a very similar way). To check that, it is enough to observe that $\psi'$ satisfies $(6.37)_1$
if and only if 
$$0 = -\psi'{}^{V_{1,1,0}}([\hat \epsilon'_3, \hat \epsilon'_1]) - 
\psi'{}^{V_{2,1,0}}([\hat \epsilon'_4, \hat \epsilon'_1]) = $$
$$ = -\psi^{V_{1,1,0}}([\hat \epsilon'_3, \hat \epsilon'_1]) - 
\psi'{}^{V_{2,1,0}}([\hat \epsilon'_4, \hat \epsilon'_1]) 
+ \sum_{j = 1}^6\left( S^{1,1,0}_j\omega^j([\hat \epsilon'_3, \hat \epsilon'_1]) +
S^{2,1,0}_j\omega^j([\hat \epsilon'_4, \hat \epsilon'_1])\right) = $$
$$ = d\psi^{V_{1,1,0}}([\hat \epsilon_3, \hat \epsilon_1]) + 
 d\psi^{V_{2,1,0}}([\hat \epsilon_4, \hat \epsilon_1])
- 2 S_1^{2,1,1} + S_4^{1,1,2} - $$
$$ - S^{1,1,0}_1d\omega^1(\hat \epsilon'_3, \hat \epsilon'_1) 
- S^{2,1,0}_1d\omega^1(\hat \epsilon'_4, \hat \epsilon'_1) = 
$$
$$ = d\psi^{V_{1,1,0}}([\hat \epsilon_3, \hat \epsilon_1]) + 
 d\psi^{V_{2,1,0}}([\hat \epsilon_4, \hat \epsilon_1])
- 2 S_1^{2,1,1} + S_4^{1,1,2}\ ,$$
where we used (6.17) and the  claim on the vanishing of certain  functions $S^{1,1,0}_j$ given 
in Lemma 6.5.\par
Now, it is clear that, even if (6.42) - (6.45) do not 
determine completely the functions $S^{1,a,1}_a$, $S^{1,a,2}_{a+2}$ and $S^{1,a,2}_{a+3}$, $a = 1,2$, 
they form a system of maximal rank, which always admits a solution.
As in the proof of the previous lemma, we may choose the functions 
$S^{A,k}_J$ so that all conditions (6.40) - (6.47) are 
satisfied at the points of a submanifold which is transversal to all fibers and then, find out what 
are the conditions in order to be able to
extend these functions to the whole trivializable open set, so that 
 (6.10) and  (3) of Lemma 6.2 are true and obtain in this way  a good 
Cartan connection. By some  straightforward computations, based on (6.12),
(6.13) and the hypothesis, we get that all this can be done implying 
 that a Cartan connection $\psi$ which satisfies 
all requirements exists and that for any other Cartan connection $\psi'$ the functions 
 $S^{A,k}_J$ satisfy the relations (6.39).
\qed
\enddemo
\medskip
\proclaim{Lemma 6.7} On any trivializable open subset $\pi^{-1}(\Cal U) \subset \P$, 
there exists a unique good Cartan connection $\psi$, which satisfies (6.17) - (6.22), (6.36) - (6.38) 
and  the following conditions for any $a = 1,2$, $j  =1,3,4$ and $k = 2,5,6$:
$$d\psi^{V_{a,1,1}}(\hat \epsilon_3, \hat \epsilon_4) = d\psi^{V_{a,2,1}}(\hat \epsilon_5, \hat \epsilon_6) = 0
\ ,\tag 6.46$$
$$d\omega^6(\hat \epsilon_j, \hat \epsilon_2) 
+ d\psi^{V_{1,2,0}}(\hat \epsilon_j, \hat \epsilon_6) - d\psi^{V_{2,2,0}}(\hat \epsilon_j, \hat \epsilon_5)
= 0\ ,\tag 6.47$$
$$d\omega^5(\hat \epsilon_j, \hat \epsilon_2) 
+ d\psi^{V_{1,2,0}}(\hat \epsilon_j, \hat \epsilon_5) + d\psi^{V_{2,2,0}}(\hat \epsilon_j, \hat \epsilon_6)
= 0\ ,\tag 6.48$$
$$d\omega^4(\hat \epsilon_k, \hat \epsilon_1) 
+ d\psi^{V_{1,1,0}}(\hat \epsilon_k, \hat \epsilon_4) - d\psi^{V_{2,1,0}}(\hat \epsilon_k, \hat \epsilon_3)
= 0\ ,\tag 6.49$$
$$d\omega^3(\hat \epsilon_k, \hat \epsilon_1) 
+ d\psi^{V_{1,1,0}}(\hat \epsilon_k, \hat \epsilon_3) + d\psi^{V_{2,1,0}}(\hat \epsilon_k, \hat \epsilon_4)
= 0\ ,\tag 6.50$$
$$d\psi^{V_{1,1,2}}(\hat \epsilon_3, \hat \epsilon_4) = 
d\psi^{V_{1,2,2}}(\hat \epsilon_5, \hat \epsilon_6) = 0\ ,
\tag 6.51$$
$$2 d\psi^{V_{2,2,0}}(\hat \epsilon_j, \hat \epsilon_2) +  d\psi^{V_{2,2,1}}(\hat \epsilon_j, \hat \epsilon_5) 
+  d\psi^{V_{1,2,1}}(\hat \epsilon_j, \hat \epsilon_6)  = 0\ ,\tag 6.52$$
$$2 d\psi^{V_{2,1,0}}(\hat \epsilon_k, \hat \epsilon_1) +  d\psi^{V_{2,1,1}}(\hat \epsilon_k, \hat \epsilon_3) 
+  d\psi^{V_{1,1,1}}(\hat \epsilon_k, \hat \epsilon_4)  = 0\ .\tag 6.53$$
\endproclaim
\demo{Proof} The proof proceeds exactly as for the previous two lemmata. Assuming that 
$\psi$ satisfies all claims of Lemma 6.6 and 6.7, let us look for a modified $\psi'$ 
which has the same  properties  and for which $(6.46)_1$ hold with $a = 1$. We get that 
the functions $S^{A,k}_i$ have to satisfy 
$$0 = d\psi'{}^{V_{1,1,1}}(\hat \epsilon'_3,\hat \epsilon'_4) = 
-  \psi^{V_{1,1,1}}([\hat \epsilon'_3,\hat \epsilon'_4]) - 
\sum_{j = 1}^6 S^{1,1,1}_jd\omega^j (\hat \epsilon'_3,\hat \epsilon'_4)  = $$
$$ = - \psi^{V_{1,1,1}}([\hat \epsilon_3,\hat \epsilon_4])  + S_3^{1,1,2} - S^{1,1,1}_1
+ \sum_{k = 5,6}ÊS^{1,1,1}_k d d\omega^k (\hat \epsilon_3,\hat \epsilon_4) = $$
$$ = 
d \psi^{V_{1,1,1}}(\hat \epsilon_3,\hat \epsilon_4) + \frac 3 2 S^{1,1,2}_3 +
\sum_{k = 5,6}ÊS^{1,1,1}_k d\omega^k (\hat \epsilon_3,\hat \epsilon_4)\ ,$$
where we used the fact that $S^{1,1,1}_j = 0$ when $j = 3,4$ and that
$d\omega^2 (\hat \epsilon_3,\hat \epsilon_4) = 0$. So we have that
$$S_3^{1,1,2} =  - \frac{2}{3} d \psi^{V_{1,1,1}}(\hat \epsilon_3,\hat \epsilon_4) +
\frac{2}{3}  \sum_{k = 5,6}ÊS^{1,1,1}_k d\omega^k (\hat \epsilon_3,\hat \epsilon_4)\ .\tag 6.54$$
On the other hand, for any $\psi'$ which is also a Cartan connection, we have to 
impose that (6.10) holds true. This implies that for any $ k = 5,6$
$$V^*_{1,1,1} (S^{1,2,2}_3) = - S^{2,1,1}_3\ ,\qquad V^*_{1,1,1} (S^{1,1,1}_k) = S^{2,1,0}_k = 0
\ ,$$
where we use the fact that $S^{2,1,0}_k = 0$ by Lemma 6.5. So, from (6.54) and
using (6.12) and (6.13), we get that any $\psi'$ which satisfies $(6.46)_1$
and (6.10) has to satisfies also
$$S^{2,1,0}_1 = S^{1,1,1}_4 =  S^{2,1,1}_3 = - V^*_{1,1,1}(S^{1,1,2}_3) = $$
$$ = - V^*_{1,1,1} \left(-
\frac{2}{3} d \psi^{V_{1,1,1}}(\hat \epsilon_3,\hat \epsilon_4) +
\frac{2}{3}  \sum_{k = 5,6}ÊS^{1,1,1}_k d d\omega^k (\hat \epsilon_3,\hat \epsilon_4)
\right) =  $$
$$ = 
\frac{2}{3} d \psi^{V_{2,1,0}}(\hat \epsilon_3,\hat \epsilon_4)
- \frac{2}{3}  \sum_{k = 5,6}ÊS^{1,1,1}_k d\omega^2 (\hat \epsilon_3,\hat \epsilon_4) = 0\ .\tag 6.55$$
With similar arguments applied to $(6.46)_2$ we get that 
$$S^{2,2,0}_2 = S^{1,2,1}_6 =  S^{2,2,1}_5 = 0 \ .\tag 6.56$$
Now, using (6.55) and (6.56), we may determine the conditions 
for which (6.47) - (6.50) are satisfied. With the same arguments of before we find that we must have
$$S^{1,2,1}_j
= - \frac{1}{3}\left(  d\omega^6(\hat \epsilon_j, \hat \epsilon_2) 
+ d\psi^{V_{1,2,0}}(\hat \epsilon_j, \hat \epsilon_6) - d\psi^{V_{2,2,0}}(\hat \epsilon_j, \hat \epsilon_5)  \right)\ ,
\tag 6.57$$ 
$$S^{2,2,1}_j
= - \frac{1}{3}\left(  d\omega^5(\hat \epsilon_j, \hat \epsilon_2) 
+ d\psi^{V_{1,2,0}}(\hat \epsilon_j, \hat \epsilon_5) + d\psi^{V_{2,2,0}}(\hat \epsilon_j, \hat \epsilon_6)  \right)\ ,
\tag 6.58$$
$$S^{1,1,1}_k
= - \frac{1}{3}\left(  d\omega^4(\hat \epsilon_k, \hat \epsilon_1) 
+ d\psi^{V_{1,1,0}}(\hat \epsilon_k, \hat \epsilon_4) - d\psi^{V_{2,1,0}}(\hat \epsilon_k, \hat \epsilon_3)  \right)
\ ,\tag 6.59$$
$$S^{2,1,1}_k
= - \frac{1}{3}\left(  d\omega^3(\hat \epsilon_k, \hat \epsilon_1) 
+ d\psi^{V_{1,1,0}}(\hat \epsilon_k, \hat \epsilon_3) + d\psi^{V_{2,1,0}}(\hat \epsilon_k, \hat \epsilon_4)  \right)
\ ,\tag 6.60$$
for any $j = 1,3,4$ and $k = 2,5,6$.
As before, a tedious but straightforward check shows that functions which satisfy 
(6.55) - (6.60) define a $\psi'$ so that
all equations (6.10) and Lemma 6.2 (3) are satisfied. So, we may replace $\psi$ with 
$\psi'$ and assume that also the conditions (6.48) - (6.50) are true. Moreover, from the previous 
discussion, we see that, if we look for another $\psi'$ which satisfies those conditions 
and so that (6.46) and (6.47) are true, we have to look among the $\gQ$-valued 1-forms
for which the functions $S^{2,1,0}_1$ and $S^{2,2,0}_2$ are vanishing, by (6.55) and 
(6.56). Moreover, by (6.57) - (6.60), we also have to suppose that 
all functions $S^{a,2,1}_j$ and $S^{a,1,1}_k$, with $j = 1, 3,4$ and $k= 2,5,6$ are 
vanishing. With these assumptions, we get  that (6.46) and (6.47) are true if and only if 
$$S_3^{1,1,2} =  - \frac{2}{3} d \psi^{V_{1,1,1}}(\hat \epsilon_3,\hat \epsilon_4)\ ,\ 
S_4^{1,1,2} =  \frac{2}{3} d \psi^{V_{2,1,1}}(\hat \epsilon_3,\hat \epsilon_4)\ , \tag 6.61$$
$$S_5^{1,2,2} = - \frac{2}{3} d \psi^{V_{1,2,1}}(\hat \epsilon_5,\hat \epsilon_6)\ ,\qquad 
S_6^{1,2,2} = \frac{2}{3} d \psi^{V_{1,2,1}}(\hat \epsilon_5,\hat \epsilon_6)\ ,\tag 6.62$$
Again, a tedious but straightforward check shows that such functions define a $\psi'$ so that
all equations (6.10) and Lemma 6.2 (3) are satisfied.
Hence we may replace $\psi$ by $\psi'$ and 
assume that now even  (6.46)  and (6.47) 
are satisfied. Moreover, any other good Cartan connection with the 
same properties is determined by functions $S^{A,k}_j$ which are all trivial 
except for the functions $S^{1,2,2}_j$ and  $S^{1,1,2}_k$, with $j = 1,3,4$ and $k = 2,5,6$, 
and the functions $S^{1,1,2}_1$ and $S^{1,2,2}_2$. \par
These  functions are uniquely determined if we require that (6.51), (6.52) and (6.53) are
true. In fact, using the vanishing of all other functions and by  the same
 arguments of before, 
we see that  this  occurs if and only if
$$S^{1,1,2}_1 = d\psi^{V_{1,1,2}}(\hat \epsilon_3, \hat \epsilon_4)\ ,\quad 
S^{1,2,2}_2 = d\psi^{V_{1,2,2}}(\hat \epsilon_5, \hat \epsilon_6)\ ,
\tag 6.63$$
$$S^{1,2,2}_j = - \frac{1}{4}
\left(
2 d\psi^{V_{2,2,0}}(\hat \epsilon_j, \hat \epsilon_2) +  d\psi^{V_{2,2,1}}(\hat \epsilon_j, \hat \epsilon_5) 
+  d\psi^{V_{1,2,1}}(\hat \epsilon_j, \hat \epsilon_6)
\right)\ ,\tag 6.64$$
$$S^{1,1,2}_k = - \frac{1}{4}
\left(
2 d\psi^{V_{2,1,0}}(\hat \epsilon_k, \hat \epsilon_1) +  d\psi^{V_{2,1,1}}(\hat \epsilon_k, \hat \epsilon_3) 
+  d\psi^{V_{1,1,1}}(\hat \epsilon_k, \hat \epsilon_4)
\right)\ ,\tag 6.65$$
for any $j = 1,3,4$ and $k = 2,5,6$.\par
Again, it can be checked that (6.10) and Lemma 6.2 (3) are satisfied and 
hence that such a $\gQ$-valued form is a Cartan connection. The uniqueness of such connection
follows immediately
from (6.63), 
(6.64) and (6.65). 
 \qed
\enddemo
The previous lemma   leads  to the 
main result of this subsection.\par
\medskip
\proclaim{Proposition 6.8} There exists a unique (globally defined) good Cartan connection $\psi_{CM}$
on $\P$ that satisfies (6.16) - (6.22), (6.25), (6.36) - (6.38), (6.46) - (6.53).  Moreover, 
the pair $(\P, \psi_{CM})$ is, generically, 
not isomorphic with the pair $(P(M), \omega_M)$ of Theorem 1.1.
\endproclaim
\demo{Proof} By the remarks after (6.10), the existence of 
a unique canonical Cartan connection which satisfies the hypothesis on 
any trivializable open set implies the existence of a globally defined Cartan connection. 
To prove the second claim, let us consider the conditions which correspond
to  (6.11). For this, 
we need to determine the basis for $\gQ^1 + \gQ^2$ which is $\B$-dual 
to the basis $\epsilon_J$, $1 \leq J\leq 6$, of $\gQ^{-1} + \gQ^{-2}$. \par
 By classical properties of graded semisimple Lie algebras (see e.g. \cite{15}, Lemma 3.15), 
the $\B$-dual element of an element in $\gQ^i$ has to be
in  $\gQ^{-i}$. In particular, the $\B$-dual element of $\epsilon_1$
is a multiple of $V_{1,1,2}$. On the other hand, 
 the 
element $V_{2,1,0}$ is  a grading element for $\gQ$, i.e. for any element $X\in \gQ^i$ in the 
 same simple subalgebra of $V_{2,1,0}$,  we have $[V_{2,1,0}, X] = i X$.
 Using this fact,  the $\B$-norm 
of $V_{2,1,0}$ is immediately computed, i.e. $\B(V_{2,1,0}, V_{2,1,0}) =  12$. So
we  also have  
$$\B(\epsilon_1, V_{1,1,2}) = \frac{1}{2}\B(\epsilon_1, [V_{2,1,0},V_{1,1,2}]) = 
- \frac{1}{2}\B([\epsilon_1,V_{1,1,2}], V_{2,1,0}) = $$
$$ = - \frac{1}{2}\B(V_{2,1,0}, V_{2,1,0}) = - 6\ .$$
This shows that the  element that is  $\B$-dual  to $\epsilon_1$ is $-\frac{1}{6} V_{1,1,2}$.
By a similar line of arguments, one can determine all the other elements $\epsilon^A$ 
of the basis which is $\B$-dual to $\epsilon_J$, $1 \leq J\leq 6$. Furthermore, 
since (6.11) is satisfied also if we rescale all vectors $\epsilon^A$ by the  
factor $-1/6$, we may assume that the vectors  $(\epsilon_A)$ are 
 $\epsilon^1 = V_{1,1,2}$, 
$\epsilon^2 = V_{1,2,2}$, $\epsilon^3 = - V_{1,1,1}$, $\epsilon^4 = V_{2,1,1}$, 
$\epsilon^5 =  - V_{1,2,1}$ and $\epsilon^6 = V_{2,2,1}$.  \par
At this point,  it is just a matter of computing and showing that 
there is one condition which is, generically, not satisfied. In fact, 
let us set $\hat \epsilon_B = \hat \epsilon_3$ 
and let us evaluate  the component along the vector $V_{1,1,2}$ of (6.11):
using (6.13), such component becomes into
$$0 = V^*_{1,1,2}(d\psi^{V_{1,1,2}}(\hat \epsilon_1, \hat \epsilon_3)) +
V^*_{2,1,1}(d\psi^{V_{1,1,2}}(\hat \epsilon_4, \hat \epsilon_3)) = $$
$$- 2 d\psi^{V_{1,1,2}}(\hat \epsilon_1, \hat \epsilon_3) - 
d\psi^{V_{1,1,2}}(\hat \epsilon_1, \hat \epsilon_3) \ .$$
By (6.46), we see that our connection satisfies (6.11) only if 
$d\psi^{V_{1,1,2}}(\hat \epsilon_1, \hat \epsilon_3) = 0$. But this 
equality is not consequences of the conditions which characterize $\psi_{CM}$, 
nor it is a consequence of the linear relations implied by the identity 
$d^2\psi_{CM} = 0$, as one can check with the help of 
a computer program like e.g. {\it Mathematica\/}. With the help of Frobenious theorem, 
one may infer that
hyperbolic manifolds, for which 
$d\psi^{V_{1,1,2}}(\hat \epsilon_1, \hat \epsilon_3) \neq 0$,  exist and 
are generic.
\qed
\enddemo
\bigskip
\subsubhead 6.2 Construction of a canonical Cartan connection on an elliptic manifold
\endsubsubhead
\medskip
Assume that $M$ is elliptic.  
In this case $\gQ = \goth{sl}_{3}(\C)$
and  the special basis $(\epsilon_i, V_{A,k})$ is given in Appendix.\par
Consider  a good Cartan connection  $\psi$ on  
a trivializable open subset $\pi^{-1}(\Cal U) \subset \P$. 
Since $\gQ$ is a complex Lie algebra, $\psi$
induces on $ \pi^{-1}(\Cal U) \subset \P$, 
 the (in general, {\it non-integrable\/}) 
complex structure, determined by the endomorphisms of tangent spaces defined by
$$J_o: T_u \P \to T_u \P\qquad \psi(J_o(X)) = i \psi(X)\ .$$
Notice that, by construction, $J_o^* \psi = i  \psi$ and that we have 
 the following relations between the 
components  of $\psi$:
$$J^*_o \psi^{\epsilon_1} = - \psi^{\epsilon_2}\ ,\quad
J^*_o\psi^{\epsilon_3} = \psi^{\epsilon_4}\ ,\quad J^*_o\psi^{\epsilon_5} = - \psi^{\epsilon_6}\ ,$$
$$J^*_o \psi^{V_{4,0}} = \psi^{V_{1,0}}\ ,\quad
J^*_o\psi^{V_{2,0}} = - \psi^{V_{3,0}}\ ,$$
$$J^*_o \psi^{V_{1,1}} = - \psi^{V_{2,1}}\ ,\quad
J^*_o\psi^{V_{3,1}} = - \psi^{V_{4,1}}\ ,J^*_o\psi^{V_{1,2}} = - \psi^{V_{2,2}}\ .$$
So, the  $\C$-valued forms on $\pi^{-1}(\Cal U) \subset \P$ defined by
$$\varpi^0 = \psi^{\epsilon_1} + i \psi^{\epsilon^2}\ ,\quad 
\varpi^1 = \psi^{\epsilon_3} - i \psi^{\epsilon^4}\ ,\quad
\varpi^2 = \psi^{\epsilon_5} + i \psi^{\epsilon^6}\ ,\tag 6.66$$
$$\Psi^{V_{I,0}} = \psi^{V_{4,0}} - i \psi^{V_{1,0}}\ ,\qquad 
\Psi^{V_{II,0}} = \psi^{V_{2,0}} + i \psi^{V_{3,0}}\ ,\qquad 
\Psi^{V_{I,1}} = \psi^{V_{1,1}} + i \psi^{V_{2,1}}\ ,\tag 6.67$$
$$
\Psi^{V_{II,1}} = \psi^{V_{3,1}} + i \psi^{V_{4,1}}\ ,\qquad
\Psi^{V_{2}} = \psi^{V_{1,2}} + i \psi^{V_{2,2}}\ ,\tag 6.68$$
determine at any tangent space $T_u\P$ a basis of complex 1-forms which are 
holomorphic  w.r.t. the complex structure $J_o: T_u\P \to T_u \P$. We may also  consider
at any point the corresponding dual holomorphic basis (i.e. a basis of elements in 
$T^{1,0} \P \subset T^\C \P$). These dual bases are given 
at any point by the 
following complex vector fields of  $T^\C \P$:
$$\hat e_0 = \frac{1}{2}\left(\hat \epsilon_1 - i \hat \epsilon_2\right) \ ,
\quad  \hat e_1 = \frac{1}{2}\left(\hat \epsilon_3 + i \hat \epsilon_4\right)\ ,
\quad \hat e_2 = \frac{1}{2}\left(\hat \epsilon_5 - i \hat \epsilon_6\right)\ ,\tag 6.69$$
$$\hat V_{I,0} = \frac{1}{2}\left( V_{4,0}^* + i V^*_{1,0}\right)\ ,\quad
\hat V_{II,0} = \frac{1}{2}\left(V_{2,0}^* - i V^*_{3,0}\right)\ ,\tag 6.70$$
$$\hat V_{I, 1} = \frac{1}{2}\left(V^*_{1,1} - i V^*_{2,1}\right)\ ,\quad 
V_{II, 1} = \frac{1}{2}\left(V_{3,1} - i V_{4,1}\right)
\ ,\quad V_2 = \frac{1}{2}\left(V_{1,2} - i V_{2,2}\right)\ .\tag 6.71$$
It is clear that any real 1-form on $\pi^{-1}(\Cal U) \subset \P$
can be expressed in terms of the real and imaginary parts of the $\C$-valued 
1-forms $\varpi^i$ and $\Psi^{V_A}$. In particular, we may express any 
good Cartan connection $\psi'$ on $\pi^{-1}(\Cal U) \subset \P$ using 
the $\varpi^i$'s and $\Psi^{V_A}$'s and some suitable complex valued functions $S^A_i$. 
Since those expressions will turn out to be very helpful for our next computations, 
we write them down, for  reader's convenience.\par
Assume that $\psi'$ is a $\gQ$-valued 1-form,
which satisfies (6.5) and i) and ii) of Definition 6.1. Then, again, we may consider the 
(non-integrable)  complex structure $J'_o$ defined by 
$i\psi'(X) = \psi(J'_o X)$ and the $\C$-valued 1-forms 
$\{ \varpi'{}^i, \Psi'{}^{V_A} \}$, which are defined as in (6.66)-(6.68) and 
are $J_o'$-holomorphic at any tangent space. We may also 
consider the associated  vector fields
$\{\hat e'_i, \hat V'_{A}\}$, which are defined as in (6.69) and (6.70). 
These vector fields and 1-forms are written in terms of  the previous one  as follows
$$\hat e_i' = \hat e_i + \sum_{A,\ell} S^{A}_i \hat V_{A} + 
\sum_{A} S^{\bar A}_i \overline{\hat V_{A}}\ ,
\qquad \hat V'_A = \hat V_A\ ,$$
$$\varpi'{}^i = \varpi^i \ ,\qquad \Psi'{}^A = \Psi^A - \sum_i S_i^A \varpi^i - 
\sum_i S_{\bar i}^A \varpi^{\bar i}\quad \left(\ \varpi^{\bar i} \= \overline{\varpi^i}\ , 
\ S_{\bar i}^A \= \overline{S_i^{\bar A}}\ 
\right)\tag 6.72$$
for some suitable $\C$-valued functions $S_i^A$ and $S_{\bar i}^A $ (which are linear
combinations of the original $\R$-valued functions $S_i^{A,\ell}$). By the  same arguments
of before, 
$\psi'$ is a Cartan connection if and only if for any vector $V_A$ amongst 
$V_{I,0} \= V_{4,0}$, $V_{II,0} \= V_{2,0}$, $V_{I,1} \= V_{1,1}$, $V_{II,1} \= V_{2,1}$ and 
$V_2 \= V_{1,2}$ and any vector 
$e_0 \= \epsilon_1$, $e_1 \= \epsilon_3$ and $e_2 \= \epsilon_5$ we have that 
$$[V_A, e_i] =  \psi'([\hat V_A, \hat e_i])\ ,\qquad 
0 = \psi'([\hat V_A, \overline{\hat e_i}])\ ,\qquad 
0 = \psi'([\overline{\hat V_A}, \hat e_i])\ .$$
By the same arguments used to prove (6.10), it follows 
that $\psi'$ is a Cartan connection if and only if 
$$\hat V_B(S^A_i) = \sum_j S^A_j \varpi^j([\hat V_B, \hat e_i]) 
- 
\sum_C S_i^C \Psi^A([\hat V_B, \hat V_C])\ ,\qquad \overline{\hat V_B}(S^A_i)  = 0 \ ,$$
$$\hat V_B(S^A_{\bar i}) = - \sum_C S_{\bar i}^C \Psi^A([\hat V_B, \hat V_C])\ ,\qquad
\overline{\hat V_B}(S^A_{\bar i}) =  \sum_j S^A_{\bar j} \overline{
\varpi^j([\hat V_B, \hat e_i])}\ .\tag 6.10'$$ 
\medskip
Finally, before going on with the next lemmata and 
construct the desired canonical Cartan connection, we want to remark that 
for any $\gQ$-valued 1-form $\psi$,
which satisfies i) and ii) of Definition 6.1 together with (6.5), 
it is possible to express $d\varpi^0 = d(\omega^1 + i \omega^2)= d(\psi^{\epsilon_1} + 
i \psi^{\epsilon_2})$ in terms of the $\C$-valued 1-forms as follows:
$$d\varpi^0 = \varpi^1 \wedge \varpi^2 + \sigma \bar \varpi^1 \wedge \bar \varpi^0
+ \tau \bar \varpi^2 \wedge \bar \varpi^0 + 2 \Psi^{V_{I,0}} \wedge \varpi^0 $$
$$
+ \ \text{linear combinations of }\ \ \left\{\ \varpi^0\wedge \bar\varpi^0\ ,
\ \varpi^0 \wedge\varpi^j\ ,\ \varpi^0 \wedge\bar \varpi^j\ ,\  1\leq j\leq 2\ \right\}\tag 6.73$$
for some suitable $\C$-valued functions $\sigma$ and $\tau$. Moreover, from the vanishing of 
$$d^2\varpi^0(\hat e_0,\overline{\hat e_1},\overline{\hat e_2})
=
d^2\varpi^0(\hat e_1,\overline{\hat e_1},\overline{\hat e_2}) =
d^2\varpi^0(\hat e_2, \overline{\hat e_1},\overline{\hat e_2}) = 0$$ 
and using (6.3), 
one can check that for any good Cartan connection $\psi$ the following identities hold:
$$d\varpi^0(\hat e_0, \overline{\hat e_0}) =
d\varpi^1(\overline{\hat e_1},\overline{\hat e_2}) = d\varpi^2(\overline{\hat e_1},\overline{\hat e_2})  = 
0\ .\tag 6.74$$
Notice that such identities could be also proved using the integrability 
of any elliptic CR structure.
\par
\bigskip
\proclaim{Lemma 6.9} On any trivializable open subset $\pi^{-1}(\Cal U) \subset \P$, 
there exists a good Cartan connection $\psi$, which 
satisfies the following conditions for any $j = 0,1,2$ together with those obtained by 
complex conjugation:
$$d\varpi^0(\hat e_0, \hat e_1) = d\varpi^0(\hat e_0, \hat e_2) = 
d\varpi^1(\hat e_1, \hat e_2) = d\varpi^2(\hat e_1, \hat e_2) = 
 0\ ,\tag 6.75$$
$$d\varpi^1(\overline{\hat e_j}, \hat e_2)
-  d\varpi^2(\overline{\hat e_j}, \hat e_1) = 
d\varpi^0(\overline{\hat e_j}, \hat e_0) +
d\varpi^1(\overline{\hat e_j}, \hat e_1) +
d\varpi^2(\overline{\hat e_j}, \hat e_2) = 0\ ,
\tag 6.76$$
$$d\varpi^1(\hat e_1, \hat e_0) = d\varpi^2(\hat e_2, \hat e_0)
= d\varpi^1(\hat e_2, \hat e_0) = d\varpi^2(\hat e_1, \hat e_0) = 0\ ,\tag 6.77$$
Moreover, if $\psi'$ is another Cartan connection with the same 
properties and expressed in terms of $\psi$ as in (6.72), then 
the corresponding 
functions $S^{I,0}_j$, $S^{II,0}_j$, with $j  = 1,2$,
and the functions $S^{I,0}_{\bar J}$, $S^{II,0}_{\bar J}$, 
with $J = 0,1,2$, 
are identically vanishing, 
while the remaining functions satisfy the relations
$$S_0^{I,0} = S_1^{II,1} = S_2^{I,1}\ ,\quad S^{II,0}_1 = - S_1^{I,1} = S_2^{II,1}\ .
\tag 6.78$$
Finally, for any good Cartan connection $\psi$, which satisfies (6.75) -(6.77) the following identities
hold for any $J = 0,1,2$:
$$ d\varpi^0(\hat e_0, \overline{\hat e_J})
= d\Psi^{I,0}(\hat e_1, \hat e_2) = d\Psi^{I,0}(\overline{\hat e_1}, 
\overline{\hat e_2})  .\tag 6.79$$
\endproclaim
\demo{Proof} As in the proof of Lemma 6.5, 
we fix a good Cartan connection $\psi$ and we consider  
a new $\gQ$-valued 1-form $\psi'$ defined by some 
functions $S^{A,k}_J$ as in (6.72). We have to show that it is possible to choose the
functions $S^{A,k}_J$ so that (6.75) and (6.76) are satisfied and so that $\psi'$ 
is a Cartan connection. By the 
same arguments as in Lemma 6.5, it turns out this occurs if and only if for $i = 1,2$
$$S_i^{I,0} =  \frac{1}{2} d\varpi^0(\hat e_0, \hat e_i)\ ,\tag 6.80$$
$$S_1^{II,0} = \frac{1}{2} d\varpi^0(\hat e_0, \hat e_2) - d\varpi^1(\hat e_1, \hat e_2)\ ,\tag 
6.81$$
$$S_2^{II,0} = - \frac{1}{2} d\varpi^0(\hat e_0, \hat e_1) - d\varpi^2(\hat e_1, \hat e_2)\ ,
\tag 6.82$$
$$S_{\bar j}^{I,0} = - \frac{1}{4}
\left(d\varpi^0(\overline{\hat e_j}, \hat e_0) + d\varpi^1(\overline{\hat e_j}, \hat e_1) +
d\varpi^2(\overline{\hat e_j}, \hat e_2)\right)\ ,\tag 6.83$$
$$ 
S_{\bar j}^{II,0} = \frac{1}{2}\left(
d\varpi^2(\overline{\hat e_j}, \hat e_1)
-  d\varpi^1(\overline{\hat e_j}, \hat e_2)
\right)
\ ,\tag 6.84$$
$$S^{I,0}_0 - S^{II,1}_1 = d\varpi^1(\hat e_1, \hat e_0)\ ,\
S^{II,0}_0 - S^{I,1}_2 = d\varpi^1(\hat e_2, \hat e_0)\ ,\tag 6.85$$
$$S^{II,0}_0 + S^{I,1}_1 = d\varpi^2(\hat e_0, \hat e_1)\ ,\
S^{I,0}_0 - S^{II,1}_2 = d\varpi^2(\hat e_2, \hat e_0)\ .\tag 6.86$$
The system given by the equations (6.80) - (6.86) has maximal rank 
and it is therefore solvable at any point. Moreover, 
with some straightforward computations based on (6.10') and 
(6.12), it can be checked that any set of functions
which satisfy that system at the points of a submanifold transversal
to all fibers can be smoothly extended to functions which 
solve all equations (6.80) - (6.86) at all points and that, at the same time,
they satisfy 
all equations (6.10') and condition (3) of Lemma 6.2. From this 
we get the existence of a Cartan connection which satisfies all 
requirements. The second claim follows immediately from 
(6.80) - (6.86). The last claim can be checked by 
evaluating 
$d^2\varpi^0(\overline{\hat e_J}, \hat e_1, \hat e_2)$, with $J = 0,1,2$, 
$d^2\varpi^0(\hat e_0, \hat e_1, \hat e_2)$ and
$d^2\varpi^0(\hat e_0, \overline{\hat e_1}, \overline{\hat e_2})$   and 
recalling that they have to be all identically vanishing 
by the triviality of the operator $d^2$.
\qed
\enddemo
\medskip
\proclaim{Lemma 6.10} On any trivializable open subset $\pi^{-1}(\Cal U) \subset \P$, 
there exists good Cartan connection $\psi$, which 
satisfies (6.75) - (6.77) and the following conditions together
with those obtained by  complex conjugation:
$$d\Psi^{V_{II,0}}(\hat e_1, \hat e_2) = 0
\ .\tag 6.87$$
$$d\Psi^{V_{II,0}}(\hat e_1, \hat e_0) + 
d\Psi^{V_{I,0}}(\hat e_2, \hat e_0) = 
d\Psi^{V_{I,0}}(\hat e_1, \hat e_0) -
d\Psi^{V_{II,0}}(\hat e_2, \hat e_0) = 0
\ .\tag 6.88$$
Moreover, if $\psi'$ is another good Cartan connections with the same 
properties and obtained from $\psi$ as in (6.72), then 
all corresponding functions $S^{A,k}_J$ are vanishing except for 
the functions which appear in $(6.78)_1$ and the 
functions $S_0^{A,1}$, $S_{\bar j}^{A,1}$ and 
$S^{I,2}_J$, with $A = I,II$, $j = 1,2$ and $J = 0,1,2$, 
among which the following relations have to be satisfied:
$$ S_0^{I,1} = - \frac{1}{2} S_1^{I,2}\ ,\qquad S_0^{II,1} = \frac{1}{2} S_2^{I,2}\ .
\tag 6.89$$
\endproclaim
\demo{Proof} Following the same line of arguments 
of the previous lemma and of Lemma 6.6 and using (6.78),
one may check that (6.87) and (6.88) are satisfied if and only if 
$$S^{II,0}_0 = S_2^{II,1} = - S^{I,1}_1 = \frac{1}{4} d\Psi^{V_{II,0}}
(\hat e_1, \hat e_2) \ ,\tag 6.90$$
$$2 S^{II,1}_0 - S^{I,2}_2 = d\Psi^{V_{II,0}}(\hat e_1, \hat e_0) + 
d\Psi^{V_{I,0}}(\hat e_2, \hat e_0)\ ,\tag 6.91$$
$$2 S^{I,1}_0 + S^{I,2}_1 = d\Psi^{V_{I,0}}(\hat e_1, \hat e_0) -
d\Psi^{V_{II,0}}(\hat e_2, \hat e_0)\ .\tag 6.92$$
As in the previous proof, we notice that 
 (6.90) - (6.92) is a maximal rank system and that, using (6.10'), 
it is possible to determine smooth solutions which solve the system 
at any point and which, at the same time, satisfy all (6.10') and 
Lemma 6.2 (3). This gives the existence of a good Cartan connection
which satisfies all requirements. The second claim  follows 
immediately from (6.90) - (6.92).
\qed
\enddemo
\medskip
\proclaim{Lemma 6.11} On any trivializable open subset $\pi^{-1}(\Cal U) \subset \P$, 
 there exists a unique  Cartan connection $\psi$, which 
satisfies (6.75) - (6.77),(6.87), (6.88)  and  the next conditions for $J = 0,1,2$, 
together with all equations that can be obtained by complex conjugation:
$$d\Psi^{V_{I,1}}(\hat e_1, \hat e_2) = d\Psi^{V_{II,1}}(\hat e_1, \hat e_2) = 0\ ,\tag 6.93$$
$$d\varpi^2(\overline{\hat e_J}, \hat e_0) + 
d\Psi^{V_{II,0}}(\overline{\hat e_J}, \hat e_2) - d\Psi^{V_{I,0}}(\overline{\hat e_J}, \hat e_1) = 0
\ ,\tag 6.94$$
$$d\varpi^1(\overline{\hat e_J}, \hat e_0) + 
d\Psi^{V_{II,0}}(\overline{\hat e_J}, \hat e_1) + d\Psi^{V_{I,0}}(\overline{\hat e_J}, \hat e_2) = 0
\ ,\tag 6.95$$
$$d\Psi^{V_{I,2}}(\hat e_1, \hat e_2) = 2 d\Psi^{V_{I,0}}(\overline{\hat e_J}, \hat e_0)
+ d\Psi^{V_{II,1}}(\overline{\hat e_J}, \hat e_1) + 
d\Psi^{V_{I,1}}(\overline{\hat e_J}, \hat e_2) = 0
\ .\tag 6.96$$
\endproclaim
\demo{Proof} Following the same line of arguments of 
the  Lemma 6.6, one can  see that any good Cartan connection for which 
(6.75) - (6.77),(6.87) and (6.88) hold, satisfies also (6.93) if and only if 
$$S^{I,2}_1 = - \frac{2}{3} d\Psi^{V_{I,1}}(\hat e_1, \hat e_2) + \frac{2}{3} 
\sum_{j = 1,2} S^{I,1}_{\bar j} d\overline{\varpi^j}(\hat e_1, \hat e_2)\ ,
\tag 6.97$$
$$S^{I,2}_2 = \frac{2}{3} d\Psi^{V_{I,1}}(\hat e_1, \hat e_2) + \frac{2}{3} 
\sum_{j = 1,2} S^{II,1}_{\bar j} d\overline{\varpi^j}(\hat e_1, \hat e_2)\ ,
\tag 6.98$$
On the other, if we assume that $\psi'$ is also a Cartan connection (and hence
the functions $S^{A,k}_J$ satisfy (6.10')), then, by (6.10), (6.12), (6.13)  and (6.78), we  must have also
that 
$$V^*_{I,1}(S^{I,1}_{\bar j}) = S^{I,0}_{\bar j} = 0\ ,$$
$$S^{I,0}_0 = S^{I,1}_2 = S^{II,1}_1 = - V^*_{I,1}(S^{I,2}_1) $$
$$ = - V^*_{I,1}\left(- \frac{2}{3} d\Psi^{V_{I,1}}(\hat e_1, \hat e_2) + \frac{2}{3} 
\sum_{j = 1,2} S^{I,1}_{\bar j} d\overline{\varpi^j}(\hat e_1, \hat e_2)\right) = 0\ .\tag 6.99$$
Using (6.99), we get that $\psi'$ satisfies also (6.94) - (6.96) if and only if 
$$S^{I,1}_{\bar J} = - \frac{1}{3}
\left(d\varpi^2(\overline{\hat e_J}, \hat e_0) + 
d\Psi^{V_{II,0}}(\overline{\hat e_J}, \hat e_2) - d\Psi^{V_{I,0}}(\overline{\hat e_J}, \hat e_1)
\right)\ ,\tag 6.100$$
$$S^{II,1}_{\bar J} = - \frac{1}{3}
\left(d\varpi^1(\overline{\hat e_J}, \hat e_0) + 
d\Psi^{V_{II,0}}(\overline{\hat e_J}, \hat e_1) + d\Psi^{V_{I,0}}(\overline{\hat e_J}, \hat e_2)
\right)\ .\tag 6.101$$
With a straightforward check of (6.10') and Lemma 6.2 (3), one can see that it is possible to 
construct a $\psi'$ which satisfies (6.99) - (6.101) and hence one for which  
(6.94) and (6.95) are true. Then, coming back to (6.93), the previous remarks
imply that the functions $S^{II,1}_{\bar j}$, $j = 1,2$, have to be vanishing and hence
that  a good Cartan connection, for which (6.93) is true,
is determined by functions $S^{A,k}_J$ such that 
$$S^{I,2}_1 = - \frac{2}{3} d\Psi^{V_{I,1}}(\hat e_1, \hat e_2) \ ,
\qquad S^{I,2}_2 = \frac{2}{3} d\Psi^{V_{I,1}}(\hat e_1, \hat e_2)\ .\tag 6.102$$
Another check based on (6.10') and Lemma 6.2 (3), shows that these conditions 
give a Cartan connection and hence that we may assume that (6.93) - (6.95)  are true. 
From (6.99) - (6.102), we also have that any other good 
Cartan connection with the same properties has all functions 
$S^{A,k}_J$  vanishing except for the  functions $S^{I,2}_{\bar J}$, $J = 0,1,2$, and for 
the function $S^{I,2}_0$. These functions 
are uniquely determined if we require that also (6.96) is satisfied. In fact 
this occurs if and only if
$$S^{I,2}_0 = d\Psi^{V_{I,0}}(\hat e_1, \hat e_2)\ ,\tag 6.103$$
$$S^{I,2}_{\bar J}Ê= - \frac{1}{4}\left(
2 d\Psi^{V_{I,0}}(\overline{\hat e_J}, \hat e_0)
+ d\Psi^{V_{II,1}}(\overline{\hat e_J}, \hat e_1) + 
d\Psi^{V_{I,1}}(\overline{\hat e_J}, \hat e_2)\right)\ . \tag 6.104$$
Another check of the validity (6.10') and of Lemma 6.2 (3) gives the existence 
of such good Cartan connection and (6.103) and (6.104) give also the uniqueness.
\qed\enddemo
\medskip
The previous lemma provides immediately the existence of 
a canonical Cartan connection for any elliptic manifold. In fact, \par
\medskip
\proclaim{Proposition 6.12} 
There exists a unique (globally defined) good Cartan connection $\psi_{CM}$
on $\P$ that satisfies (6.74), (6.75) - (6.77), (6.79), (6.87), (6.88), (6.94) - (6.96). 
 Moreover, the pair $(\P, \psi_{CM})$ is, generically, 
not isomorphic  with the 
pair $(P(M), \omega_M)$ of Theorem 1.1.
\endproclaim
\demo{Proof} By the remarks after (6.10), the existence of 
a unique canonical Cartan connection which satisfies the hypothesis on 
any trivializable open set implies the existence of a globally defined Cartan connection.
For the second claim, it suffices to observe that the 
component of (6.11) along the vector $V_{I,2}$ with $\epsilon_B = \epsilon_3$
is generically not satisfied. Since the computations to explicitate 
such component and to check that the corresponding condition
is not satisfied are are very similar to those 
given in the proof of Proposition 6.8, we omit them.
\qed
\enddemo
\bigskip
\bigskip

\head Appendix
\endhead
\centerline{\it A special basis for $\gQ$ (hyperbolic case)}
\medskip
It is known (see e.g. \cite{16}) that in this case $\gQ$ is a Lie algebra isomorphic to 
$$\goth{su}_{2,1} \oplus \goth{su}_{2,1}  
 \simeq \left\{ \left( \matrix 
\smallmatrix b_1 - i\frac{b'_1}{2} & \alpha_1  & a_1 \\
2 i \bar \beta_1 & i b'_1  & 2 i \bar \alpha_1 \\
c_1 & \beta_1 & - b_1 - i\frac{b'_1}{2} \endsmallmatrix & 0 \\
0 & 
\smallmatrix b_2 - i\frac{b'_2}{2} & \alpha_2  & a_2 \\
2 i \bar \beta_2 & i b'_2  & 2 i \bar \alpha_2 \\
c_2 & \beta_2 & - b_2 - i\frac{b'_2}{2} \endsmallmatrix\endmatrix \right)
\matrix \text{with}\ 
a_i, b_i, b'_i,  c_i \in \R\\
\text{and}\ \alpha_i, \beta_i \in \C\ \endmatrix\right\}$$
The basis we consider for such 
Lie algebra is the following:
$$\epsilon_1 = \left( \matrix 
\smallmatrix 0 & 0  & - 1 \\
0 & 0  & 0 \\
0 & 0 & 0 \endsmallmatrix & 0 \\
0 & 
\smallmatrix 0 & 0  & 0  \\
0 & 0 & 0 \\
0 & 0  & 0   \endsmallmatrix \endmatrix \right) \ ,\qquad \epsilon_2 = 
\left( \matrix 
\smallmatrix 0 & 0  & 0  \\
0 & 0 & 0 \\
0 & 0  & 0   \endsmallmatrix
& 0 \\
0 & 
\smallmatrix 0 & 0  & - 1\\
0 & 0  & 0 \\
0 & 0 & 0   \endsmallmatrix \endmatrix \right)\ ,$$
$$\epsilon_3 = \left( \matrix 
\smallmatrix 0 & \frac{1}{2}  &  0 \\
0 & 0  & i\\
0 & 0 & 0   \endsmallmatrix & 0 \\
0 & 
\smallmatrix 0 & 0  & 0  \\
0 & 0 & 0 \\
0 & 0  & 0 \endsmallmatrix \endmatrix \right) \ ,\qquad \epsilon_5 = 
\left( \matrix 
\smallmatrix 0 & 0  & 0\\
0 & 0 & 0 \\
0 & 0  & 0   \endsmallmatrix
& 0 \\
0 & 
\smallmatrix 0 & \frac{1}{2}  &  0 \\
0 & 0  & i\\
0 & 0 & 0   \endsmallmatrix \endmatrix \right)\ ,$$
$$\epsilon_4 = \left( \matrix 
\smallmatrix 0 & \frac{i}{2}   & 0 \\
0 & 0  & 1\\
0 & 0 & 0  \endsmallmatrix & 0 \\
0 & 
\smallmatrix 0 & 0  & 0 \\
0 & 0 & 0 \\
0 & 0 & 0 \endsmallmatrix \endmatrix \right) \ ,\qquad \epsilon_6 = 
\left( \matrix 
\smallmatrix 0 & 0  & 0 \\
0 & 0 & 0 \\
0 & 0  & 0  \endsmallmatrix
& 0 \\
0 & 
\smallmatrix 0 & \frac{i}{2} & 0 \\
0 & 0  & 1\\
0 & 0 & 0  \endsmallmatrix \endmatrix \right)\ ,$$
$$V_{1,1,0} = \left( \matrix 
\smallmatrix \frac{i}{3} & 0  & 0  \\
0 & - \frac{2}{3} i  & 0 \\
0 & 0  &  \frac{i}{3} 
 \endsmallmatrix & 0 \\
0 & 
\smallmatrix 0 & 0  & 0 \\
0 & 0 & 0 \\
0 & 0  & 0   \endsmallmatrix \endmatrix \right) \ ,\qquad V_{1,2,0} = 
\left( \matrix 
\smallmatrix 0 & 0  & 0 \\
0 & 0 & 0 \\
0 & 0  & 0  \endsmallmatrix
& 0 \\
0 & 
\smallmatrix  \frac{i}{3} & 0  & 0  \\
0 &  -\frac{2}{3} i  & 0 \\
0 & 0  &  \frac{i}{3}  
 \endsmallmatrix \endmatrix \right)\ ,$$
$$V_{2,1, 0} = \left( \matrix 
\smallmatrix -1 & 0  & 0  \\
0 & 0  & 0 \\
0 & 0 & 1  \endsmallmatrix & 0 \\
0 & 
\smallmatrix 0 & 0  & 0  \\
0 & 0 & 0 \\
0 & 0  & 0   \endsmallmatrix \endmatrix \right) \ ,\qquad V_{2,2, 0} = 
\left( \matrix 
\smallmatrix 0 & 0  & 0  \\
0 & 0 & 0 \\
0 & 0  & 0  \endsmallmatrix
& 0 \\
0 & 
\smallmatrix -1 & 0  & 0 \\
0 & 0 & 0 \\
0 & 0 & 1  \endsmallmatrix \endmatrix \right)\ ,$$
$$V_{1, 1, 1} = \left( \matrix 
\smallmatrix 0 & 0  & 0\\
1 & 0   & 0 \\
0 & - \frac{i}{2} & 0  \endsmallmatrix & 0 \\
0 & 
\smallmatrix 0 & 0  & 0  \\
0 & 0 & 0 \\
0 & 0  & 0  \endsmallmatrix  \endmatrix \right) \ ,\qquad V_{1, 2, 1} = 
\left( \matrix
\smallmatrix
0 & 0 & 0 \\
0 & 0 & 0 \\
0 & 0  & 0   \endsmallmatrix & 0 \\
0 & 
\smallmatrix 0 & 0  & 0  \\
1 & 0 & 0  \\
0 & - \frac{i}{2} & 0  \endsmallmatrix 
\endmatrix \right)\ ,$$
$$V_{2, 1, 1} = \left( \matrix 
\smallmatrix 0 & 0  & 0  \\
i & 0   & 0  \\
0 & -\frac{1}{2}  &  0  \endsmallmatrix & 0 \\
0 & 
\smallmatrix 0 & 0  & 0  \\
0 & 0 & 0 \\
0 & 0  & 0   \endsmallmatrix \endmatrix \right) \ ,\qquad V_{2, 2, 1} = 
\left( \matrix 
\smallmatrix 0 & 0  & 0  \\
0 & 0 & 0 \\
0 & 0  & 0   \endsmallmatrix & 0 \\
0 & 
\smallmatrix 0 & 0  & 0  \\
i & 0   & 0   \\
0 & - \frac{1}{2}   & 0  \endsmallmatrix \endmatrix \right)\ ,$$
$$V_{1, 1, 2} = \left( \matrix 
\smallmatrix 0 & 0  & 0 \\
 0  & 0  & 0  \\
1 & 0 & 0   \endsmallmatrix & 0 \\
0 & 
\smallmatrix 0 & 0  & 0  \\
0 & 0 & 0 \\
0 & 0  & 0   \endsmallmatrix \endmatrix \right) \ ,\qquad V_{1, 2, 2} = 
\left( \matrix 
\smallmatrix 0 & 0  & 0  \\
0 & 0 & 0 \\
0 & 0  & 0   \endsmallmatrix & 0 \\
0 & 
\smallmatrix 0 & 0  & 0 \\
 0  & 0  & 0  \\
1 & 0 & 0   \endsmallmatrix \endmatrix \right)\ ,$$
\bigskip
\centerline{\it A special basis for $\gQ$  (elliptic case)} 
\medskip
It is known (see e.g. \cite{16}) that in this case $\gQ$ is a Lie algebra isomorphic to 
$$\goth{sl}_{3}(\C)  =  \left\{ \left(\matrix
 \alpha_{11} &  \alpha_{12} &  \alpha_{13} \\
\alpha_{21} &  \alpha_{22} &  \alpha_{23} \\
\alpha_{31} &  \alpha_{32} &  \alpha_{33} \\
\endmatrix \right)\ ,\ 
\ \alpha_{ij} \in \C \ ,\ \alpha_{11} + 
\alpha_{22} + \alpha_{33}  = 0\ \right\}\ .$$
The basis we consider for such 
Lie algebra is the following:
$$\epsilon_1 = \left( \matrix 
0 & 0  &  1 \\
0 & 0  & 0 \\
0 & 0 & 0  \endmatrix \right) \ ,\qquad \epsilon_2 = 
\left( \matrix 
0 & 0  & i \\
0 & 0  & 0 \\
0 & 0 & 0  \endmatrix \right)\ ,$$
$$\epsilon_3 = \left( \matrix 
0 & \frac{1}{2} & 0 \\
0 & 0  & i\\
0 & 0 & 0  \endmatrix \right) \ ,\qquad  \epsilon_5 = \left( \matrix 
0 & - \frac{i}{2} & 0 \\
0 & 0  &  - 1\\
0 & 0 & 0  \endmatrix \right) \ ,$$
$$\epsilon_4 = 
\left( \matrix 
0 &  - \frac{i}{2} & 0 \\
0 & 0  &  1\\
0 & 0 & 0  \endmatrix \right) \ ,\qquad  \epsilon_6 = 
\left( \matrix 
0 &  \frac{1}{2} & 0 \\
0 & 0  &  - i \\
0 & 0 & 0  \endmatrix \right) \ ,$$
$$V_{1,0} = 
\left( \matrix 
i & 0 & 0 \\
0 & 0  & 0\\
0 & 0 & - i \endmatrix \right)
\ ,\ V_{2,0} =
\left( \matrix 
- \frac{i}{3} & 0 & 0 \\
0 & \frac{2 i}{3}  & 0\\
0 & 0 & - \frac{i}{3}  
\endmatrix \right)
\ ,$$
$$ V_{3,0} = 
\left( \matrix 
\frac{1}{3} & 0 & 0 \\
0 & - \frac{2}{3}  & 0\\
0 & 0 & \frac{1}{3}  \endmatrix \right)
,\ 
V_{4,0} = 
\left( \matrix 
 - 1 & 0 & 0 \\
0 &  0 & 0\\
0 & 0 &  1 \endmatrix \right) 
,$$
$$V_{1, 1} = \left( \matrix 
0 & 0 & 0 \\
1 & 0  & 0\\
0 &  - \frac{i}{2} & 
 0  \endmatrix \right) 
\ ,\qquad V_{3, 1} = \left( \matrix 
0 & 0 & 0 \\
 - i & 0  & 0\\
0 &  \frac{1}{2} & 
 0  \endmatrix \right) \ ,$$
$$V_{2, 1} = \left( \matrix 
0 & 0 & 0 \\
 i & 0  & 0\\
0 & \frac{1}{2} & 
 0  \endmatrix \right) \ ,\qquad V_{4, 1} = 
\left( \matrix 
0 & 0 & 0 \\
1 & 0  & 0\\
0 & \frac{i}{2} & 
 0  \endmatrix \right)\ ,$$
$$V_{1,2} = \left( \matrix 
0 & 0 & 0 \\
0 & 0  & 0\\
-1 & 
 0  &  0  \endmatrix \right)  \ ,\qquad  
V_{2, 2} = \left( \matrix 
0 & 0 & 0 \\
0 & 0  & 0\\
 - i & 
 0  &  0  \endmatrix \right) $$

\bigskip
\bigskip
\font\smallsmc = cmcsc8
\font\smalltt = cmtt8
\font\smallit = cmti8
\hbox{\parindent=0pt\parskip=0pt
\vbox{\baselineskip 9.5 pt \hsize=3.1truein
\obeylines
{\smallsmc
Gerd Schmalz 
Mathematisches Institut
Rheinische Friederich-Wilhelms-Universit\"at
Beringstrasse 1
53115 Bonn
GERMANY
}\medskip
{\smallit E-mail address}\/: {\smalltt  schmalz\@math.uni-bonn.de
}
}
\hskip 0.0truecm
\vbox{\baselineskip 9.5 pt  \hsize=3.7truein
\obeylines
{\smallsmc
Andrea Spiro
Dipartimento di Matematica e Informatica
Universit\`a di Camerino
Via Madonna delle Carceri
I-62032 Camerino (Macerata)
ITALY
}\medskip
{\smallit E-mail address}\/: {\smalltt andrea.spiro\@unicam.it}
}
}
\enddocument

\bye